\documentclass[12pt,a4paper]{article}
\usepackage{amsmath,amsthm,amsfonts,graphicx,psfrag,amssymb,dsfont,mathrsfs}
\usepackage{minitoc,color}

\usepackage[colorlinks]{hyperref}
\usepackage[alphabetic,initials]{amsrefs}
\usepackage{bbm}
\usepackage{dsfont}
\usepackage[ngerman,english]{babel}
\usepackage{stackengine}
\usepackage{tikz}
\usepackage{enumitem}

\usetikzlibrary{matrix,arrows}

\newcommand{\crit}{\operatorname{crit}}
\newcommand{\sign}{\operatorname{sign}}

\def\br#1{\left\{#1\right\}}

\theoremstyle{plain}
\newtheorem{thm}{Theorem}[section]

\newtheorem{prop}[thm]{Proposition}
\newtheorem{cor}[thm]{Corollary}
\newtheorem{lemma}[thm]{Lemma}

\theoremstyle{definition}

\newtheorem{notation}[thm]{Notation}
\newtheorem{defn}[thm]{Definition}

\newtheorem{disclaimer}[thm]{Disclaimer}
\newtheorem{remark}[thm]{Remark}
\newtheorem{conj}[thm]{Conjecture}

\usepackage{scalerel}

\newcommand\reallywidehat[1]{\arraycolsep=0pt\relax%
\begin{array}{c}
\stretchto{
  \scaleto{
    \scalerel*[\widthof{\ensuremath{#1}}]{\kern-.5pt\bigwedge\kern-.5pt}
    {\rule[-\textheight/2]{1ex}{\textheight}} 
  }{\textheight} %
}{0.8ex}\\           
#1\\                 
\rule{-1ex}{0ex}
\end{array}
}

\date{}
\title{SFT computations and intersection theory in higher-dimensional contact manifolds}
\author{Agustin Moreno}

\begin{document}

\maketitle
\thispagestyle{empty}

\begin{abstract}
We construct infinitely many non-diffeomorphic examples of $5$-dimensional contact manifolds which are tight, admit no strong fillings, and do not have Giroux torsion. We obtain obstruction results for symplectic cobordisms, for which we give a proof not relying on the polyfold abstract perturbation scheme for SFT. These results are part of the author's PhD thesis \cite{Mo2}, and are the first applications of higher-dimensional Siefring intersection theory for holomorphic curves and hypersurfaces, as outlined in \cite{Mo2, MS}, as a prequel of \cite{Sie2} . \end{abstract}

\addtocontents{toc}{\protect\sloppy}
\tableofcontents

\section{Introduction}

This paper is a followup to \cite{Mo}, where we address the general problem of constructing ``interesting'' families of contact structures in higher dimensions, together with developing general computational techniques for SFT-type invariants. 

In \cite{Mo} (or \cite{Mo2}), we constructed families of contact manifolds $M$ in any odd-dimension, which have non-zero and finite algebraic torsion, in the sense of Latschev-Wendl \cite{LW}. In particular, they are tight and do not admit strong symplectic fillings. We also established, in higher dimensions, that Giroux torsion (in the sense of \cite{MNW}) implies algebraic $1$-torsion. Our examples present a geometric structure which we call a \emph{spinal open book decomposition} or SOBD (based on the 3-dimensional version of \cite{LVHMW}). This geometric structure, which one could call \emph{partially planar}, ``supports'' a suitable contact structure, and is of a certain type which can be ``detected'' algebraically by the SFT machinery.

In this paper, we discuss the same examples of \cite{Mo}, from a ``dual'' point of view. While keeping the same underlying geometric decomposition, the pieces of the decomposition have reversed roles. The resulting supported contact structure is isotopic to the original one, which is the expected behaviour, in terms of the expected ``Giroux''-type correspondence in the setting of SOBDs. However, the associated holomorphic data, originally a ($2$-dimensional) finite energy foliation, is replaced by a ($2$-codimensional) foliation by holomorphic hypersurfaces. While the holomorphic curve invariants should depend only on the isotopy class of the supported contact structure, the two points of view complement each other, from a computational point of view. For instance, in \cite{Mo}, we used the finite energy foliation to bound the order of algebraic torsion from above. In this paper, we use the dual foliation to argue that the bounds from \cite{Mo} are not, in general, necessarily the optimal. 

After carrying out a detailed construction of our examples in any odd dimension, we focus on a family of $5$-dimensional particular cases, in order to illustrate the techniques developed for the general case. For these, we study their SFT differential, up to order $2$. While our original aim was to show that the order of algebraic torsion is strictly greater than $1$, the result turned out to be quite unexpected. After classifying all possible configurations of holomorphic buildings which might contribute to algebraic $1$-torsion (there are precisely $35$ of them), all of them come in cancelling pairs, except only \emph{one} of them. While we cannot rigorously prove that counting elements in this configuration, which we call \emph{sporadic}, gives algebraic $1$-torsion, we give a heuristic argument based on string topology \cite{CL09} showing that it very likely does.   

However, even though our computation is not the expected one, what is remarkable is that all the expected corollaries \emph{still hold}. From our knowledge of the SFT differential of our $5$-dimensional examples, we show that they do not have Giroux torsion (which would follow if they didn't have algebraic $1$-torsion). We also obtain non-existence results of symplectic cobordisms for these $5$-dimensional contact manifolds. 

The failure of our computation can be given the following conjectural interpretation. Out of the $35$ aforementioned configurations, we have two types: cylinders, or tori with one puncture. The former all come in cancelling pairs, and among the latter, only the sporadic one does not. For the corollaries, all we use is the cancellation phenomenom for \emph{cylinders}, whereas the tori play no role whatsoever. In contrast, both are taken into consideration for algebraic $1$-torsion. This suggests the existence of an invariant more subtle than algebraic torsion, which we suspect would be obtained from the \emph{rational} SFT, rather than the whole SFT, and yet needs to be discovered. Morally, in the proof of the corollaries we would be exploiting the properties of this hypothetical invariant.

A significant difficulty, in practice, is the lack of a higher-dimensional intersection theory between punctured holomorphic curves, in the sense of Richard Siefring. In dimensions 3 and 4, this is a good tool for proving that holomorphic curves with certain prescribed asymptotics are unique, and therefore one knows exactly what to count. In higher dimensions, although it doesn't make sense to count intersections between holomorphic curves, it does make sense to count intersections between holomorphic curves and hypersurfaces. But, as in the $4$-dimensional situation, neither curves or hypersurfaces are closed, so that intersections coming ``from infinity'' need to be considered.

In this paper, we will obtain the first applications of the basic intersection theory between curves and hypersurfaces which are asymptotically cylindrical in a well-defined sense. This is outlined in Appendix C in \cite{Mo2}, co-written with Richard Siefring (to appear as an independent article \cite{MS}) and is a prequel of his upcoming work \cite{Sie2}, generalizing his results of \cite{Sie11}. We shall use the results of \cite{MS, Mo2} to restrict the behaviour of holomorphic curves in our examples. They necessarily lie in the leaves of a suitable codimension-$2$ holomorphic foliation. This will be crucial for extracting information from the SFT of our examples.

Moreover, we need to understand the number of ways certain holomorphic building configurations may glue to honest curves, which we intend to count, after making $J$ generic. For this, we make use of obstruction bundles in the sense of Hutchings--Taubes \cite{HT1,HT2}, a fairly non-trivial gadget to deal with in practice. The idea is to count the number of honest curves obtained by gluing buildings, by algebraically counting the zero set of a section of an obstruction bundle. While we will not explicitly compute those numbers, and just prove existence results in suitable cases, the symmetries in the setup imply that there are cancellations which can be exploited to obtain our results. Under the assumption that obstruction bundles exist leaf-wise, we show their existence in our examples by careful analytical considerations, exploiting the fact that curves lies in hypersurfaces of the foliation.

\subparagraph*{On the invariant.} The invariant we will use, algebraic torsion, was defined in \cite{LW}, and is a contact invariant taking values in $\mathbb{Z}^{\geq 0}\cup \{\infty\}$. It was introduced, using the machinery of \emph{Symplectic Field Theory}, as a quantitative way of measuring non-fillability, giving rise to a ``hierarchy of fillability obstructions'', cf.\ \cite{Wen2}. At least morally, $0$-torsion should correspond to overtwistedness, whereas $1$-torsion is implied by Giroux torsion (the converse is not true). Having $0$-torsion is actually equivalent to being \emph{algebraically overtwisted}, which means that the contact homology, or equivalently its SFT, vanishes (Proposition 2.9 in \cite{LW}). This is well-known to be implied by overtwistedness, but the converse is still wide open.

The key fact about this invariant is that it behaves well under exact symplectic cobordisms, which implies that the concave end inherits any order of algebraic torsion that the convex end has. Thus, algebraic torsion may be also thought of as an obstruction to the existence of exact symplectic cobordisms. In particular, it serves as an obstruction to symplectic fillability. Moreover, there are connections to dynamics: any contact manifold with finite torsion satisfies the Weinstein conjecture (i.e.\ there exist closed Reeb orbits for every contact form).  

One should mention that there are other notions of algebraic torsion in the literature which do not use SFT, but which are only 3-dimensional (see \cite{KMvhMW} for the version using Heegard Floer homology, or the appendix in \cite{LW} by Hutchings, using ECH).

\subparagraph*{Statement of results.}

For the SFT setup, we follow \cite{LW}, where we refer the reader for more details. We will take the SFT of a contact manifold $(M,\xi)$ (with coefficients) to be the homology $H^{SFT}_*(M, \xi; \mathcal{R})$ of a $\mathbb{Z}_2$-graded unital $BV_\infty$-algebra $(\mathcal{A}[[\hbar]], \mathbf{D}_{SFT})$ over the group ring $R_{\mathcal{R}}:=\mathbb{R}[H_2(M;\mathbb{R})/\mathcal{R}]$, for some linear subspace $\mathcal{R}\subseteq H_2(M;\mathbb{R})$. Here, $\mathcal{A}=\mathcal{A}(\lambda)$ has generators $q_\gamma$ for each good closed Reeb orbit $\gamma$ with respect to some nondegenerate contact form $\lambda$ for $\xi$, $\hbar$ is an even variable, and the operator $$\mathbf{D}_{SFT}: \mathcal{A}[[\hbar]] \rightarrow\mathcal{A}[[\hbar]]$$ is defined by counting rigid solutions to a suitable abstract perturbation of a $J$-holomorphic curve equation in the symplectization of $(M, \xi)$. It satisfies

\begin{itemize}
 \item $\mathbf{D}_{SFT}$ is odd and squares to zero,
 \item $\mathbf{D}_{SFT}(1) = 0$, and
 \item $\mathbf{D}_{SFT} = \sum_{k \geq 1} D_k\hbar^{k-1},$
\end{itemize}

where $D_k : \mathcal{A} \rightarrow \mathcal{A}$ is a differential operator of order $\leq k$, given by
$$
D_k = \sum_{\substack{\Gamma^+,\Gamma^-, g,d\\ |\Gamma^+|+g=k}}\frac{n_g(\Gamma^+,\Gamma^-,d)}{C(\Gamma^-,\Gamma^+)}q_{\gamma_1^-}\dots q_{\gamma_{s^-}^-}z^d\frac{\partial}{\partial q_{\gamma_1^+}}\dots \frac{\partial}{\partial q_{\gamma_{s^+}^+}}
$$
The sum ranges over all non-negative integers $g \geq 0$, homology classes $d \in H_2(M; \mathbb{R})/\mathcal{R}$ and ordered (possibly empty) collections of good closed Reeb orbits $\Gamma^\pm = (\gamma_1^{\pm},\dots,\gamma_{s^\pm}^\pm)$ such that $s^+ + g = k$. After a choice of spanning surfaces as in \cite{EGH} (p. 566, see also p. 651), the projection to $M$ of each finite
energy holomorphic curve $u$ can be capped off to a 2-cycle in $M$, and so it gives rise to a homology class $[u]\in H_2(M)$, which we project to define $\overline{[u]} \in H_2(M; \mathbb{R})/\mathcal{R}$. The number $n_g(\Gamma^+,\Gamma^-,d) \in \mathbb{Q}$ denotes the count of (suitably perturbed) holomorphic curves of genus $g$ with positive asymptotics $\Gamma^+$ and negative asymptotics $\Gamma^-$ in the homology class $d$, including asymptotic markers as explained in \cite{EGH}, or \cite{Wen3}, and including rational weights arising from automorphisms. $C(\Gamma^-, \Gamma^+) \in \mathbb{N}$ is a combinatorial factor defined as $C(\Gamma^-, \Gamma^+) = s^-!s^+!\kappa_{\gamma_1^-}\dots\kappa_{\gamma_{s^-}^-}$, where $\kappa_\gamma$ denotes the covering multiplicity of the Reeb orbit $\gamma$.

The most important special cases for our choice of linear subspace $\mathcal{R}$ are $\mathcal{R} = H_2(M; \mathbb{R})$ and $\mathcal{R} = \{0\}$, called the \emph{untwisted} and \emph{fully twisted} cases respectively, and $\mathcal{R} = \ker \Omega$ with $\Omega$ a closed 2-form on $M$. We shall abbreviate the latter case as $H^{SFT}_*(M, \xi;\Omega) := H^{SFT}_*(M, \xi; \ker \Omega)$, and the untwisted case simply by $H^{SFT}_*(M, \xi):=H^{SFT}_*(M, \xi;H_2(M;\mathbb{R}))$.

\begin{defn} Let $(M, \xi)$ be a closed manifold of dimension $2n+1$ with a positive, co-oriented contact structure. For any integer $k \geq 0$, we say that $(M, \xi)$ has $\Omega$-twisted algebraic torsion of order $k$ (or $\Omega$-twisted $k$-torsion) if $[\hbar^k] = 0$ in $H^{SFT}_*(M, \xi;\Omega)$. If this is true for all $\Omega$, or equivalently, if $[\hbar^k] = 0$ in $H^{SFT}_*(M, \xi;\{0\})$, then we say that $(M, \xi)$ has fully twisted algebraic $k$-torsion. 
\end{defn}

We will refer to \emph{untwisted} $k$-torsion to the case $\Omega=0$, in which case $R_\mathcal{R}=\mathbb{R}$ and we do not keep track of homology classes. Whenever we refer to torsion without mention to coefficients we will mean the untwisted version. We will say that, if a contact manifold has algebraic $0$-torsion for every choice of coefficient ring, then it is \emph{algebraically overtwisted}, which is equivalent to the vanishing of the SFT, or its contact homology. By definition, $k$-torsion implies $(k+1)$-torsion, so we may define its \emph{algebraic torsion} to be
$$
AT(M,\xi;\mathcal{R}):=\min\{k\geq 0:[\hbar^k]=0\}\in \mathbb{Z}^{\geq 0}\cup\{\infty\},
$$
where we set $\min \emptyset=\infty$. We denote it by $AT(M,\xi)$, in the untwisted case.

\vspace{0.5cm}

Examples of 3-dimensional contact manifolds with any given order of torsion $k-1$, but not $k-2$, were constructed in \cite{LW}. We will consider a generalization of their examples. Consider $\Sigma$ a surface of genus $g$ which is divided into two pieces $\Sigma_+$ and $\Sigma_-$ along some dividing set of simple closed curves $\Gamma$ of cardinality $k$, where the latter has genus $0$, and the former has genus $g-k+1$. Consider also a \emph{closed} $(2n-1)$-manifold $Y$, such that $Y\times [-1,1]$ admits the structure of a Liouville domain (which we call a \emph{cylindrical Liouville semi-filling}). Let $M_g:=Y\times \Sigma$.  

We will fix the following notation:

\begin{notation}\label{not} Throughout this document, the symbol $I$ will be reserved for the interval $[-1,1]$.
\end{notation}

\begin{figure}[t]\centering \includegraphics[width=0.70\linewidth]{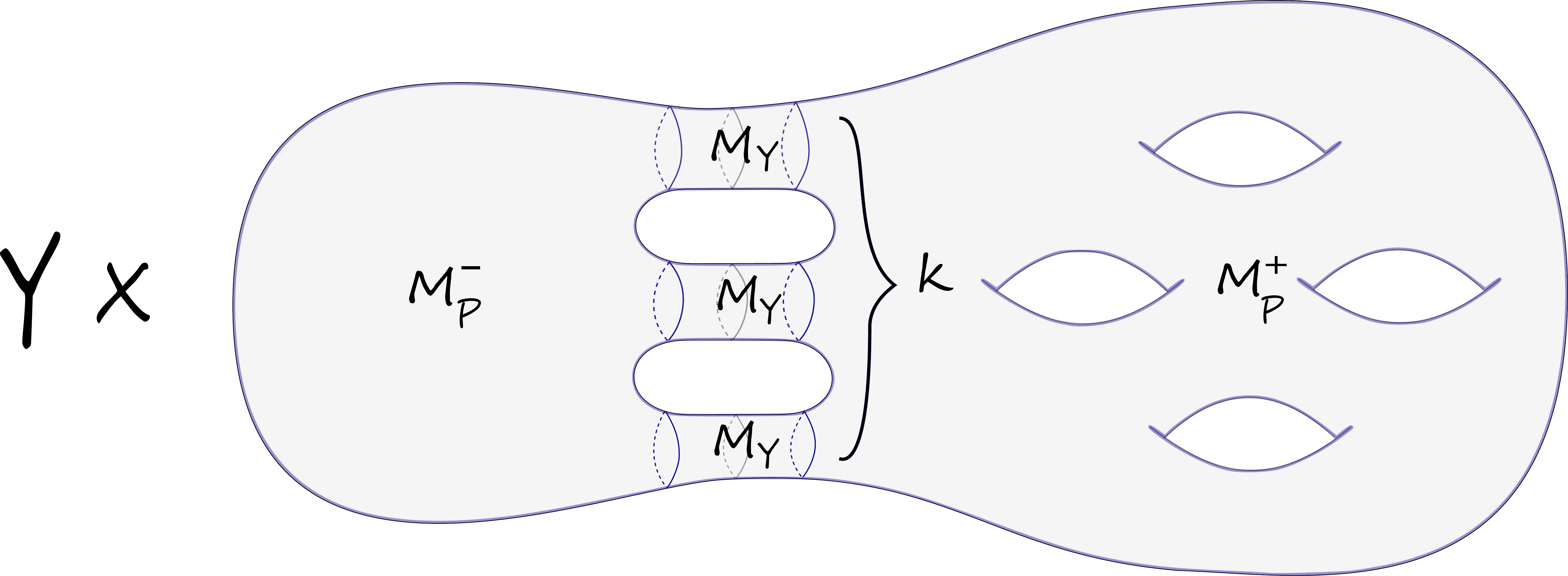}
\caption{\label{themanifold} The SOBD structure in $M$.}
\end{figure}

We can adapt the construction of the contact structures in \cite{LW} to our models. We decompose the manifold $	M=M_g=Y\times \Sigma$ into three pieces
$$
M_g=M_Y \bigcup M_P^\pm,
$$
where $M_Y = \bigsqcup^k Y \times I \times S^1$ is the \emph{spine}, and $M_P^\pm = Y \times \Sigma_\pm$ is the \emph{paper} (see Figure \ref{themanifold}). 
We have natural fibrations $$\pi_Y: M_Y \rightarrow Y \times I$$$$ \pi_P^\pm: M_P^{\pm} \rightarrow Y_\pm,$$ with fibers $S^1$ and $\Sigma_\pm$, respectively, and they are compatible in the sense that
$$
\partial((\pi_P^\pm)^{-1}(pt))=\bigsqcup^k \pi_Y^{-1}(pt)
$$

While $\pi_Y$ has a Liouville domain as base, and a contact manifold as fiber, the situation is reversed for $\pi_P^\pm$, which has contact base, and Liouville fibers. This is a prototypical example of a SOBD. While we will not include a formal definition of this notion, the reader is invited to consult \cite{Mo2} for a tentative one.

Using this decomposition, we can construct a contact structure $\xi_k$ which is a small perturbation of the ``confoliation''-type stable Hamiltonian structure $\xi_\pm \oplus T\Sigma_\pm$ along $M_P^\pm$, and is a \emph{contactization} for the Liouville domain $(Y\times I,\epsilon d\alpha)$ along $M_Y$, for some small $\epsilon>0$. This means that it coincides with $\ker(\epsilon\alpha + d\theta)$, where $\theta$ is the $S^1$-coordinate. This was done in detail in \cite{Mo}, where the following was proved: 

\begin{thm}\cite{Mo}\label{thm5d}
For any $k\geq 1$, and $g\geq k$, the $(2n+1)$-dimensional contact manifolds $(M_g=Y\times \Sigma,\xi_k)$ satisfy $AT(M_g,\xi_k)\leq k-1$. In particular, it  is not (strongly) symplectically fillable. Moreover, if $(Y,\alpha_\pm)$ are hypertight, and $k\geq 2$, the corresponding contact manifold $(M_g,\xi_k)$ is also hypertight. In particular, $AT(M_g,\xi_k)>0$, and it is tight.  
\end{thm}

In this paper, we consider the ``dual'' SOBD, where the roles of the fibrations are reversed:
$$\pi^*_Y: M_Y \rightarrow S^1$$
$$ (\pi_P^\pm)^*: M_P^{\pm} \rightarrow \Sigma_\pm,$$
so that $\pi^*_Y$ is a Liouville fibration over a contact base, and $(\pi_P^\pm)^*$, a contact fibration over a Liouville base. Observe that we can do this due to the absence of monodromy (the fibrations are trivial).

The associated contact structure $\xi_k^*$ is now a perturbation of the (integrable) stable Hamiltonian structure $T(Y \times I)$ along $M_Y$, which is now the \emph{paper}, and isotopic to a ``$(Y_\pm,\alpha_\pm)$-contactization'' of the Liouville domain $(\Sigma_\pm,\epsilon\lambda_\pm)$ along $M_P^\pm$, which is now the spine (for some Liouville form $\lambda_\pm$). More concretely, $\xi_k^*$ is isotopic, along $M_P^\pm$, to $\ker(\epsilon \lambda_\pm + \alpha_\pm)$. One can check explicitly that $\xi_k$ and $\xi_k^*$ are isotopic contact structures \cite{Mo2}. We will drop the $*$ from the notation, since our statements only depend on the isotopy type of the contact structure. 

\vspace{0.5cm}

The authors of \cite{MNW} define a generalized higher-dimensional version of the notion of Giroux torsion. This notion is defined as follows: consider $(Y,\alpha_+,\alpha_-)$ a Liouville pair on a closed manifold $Y^{2n-1}$, i.e.\ the $1$-form $\beta=\frac{1}{2}(e^s\alpha_++e^{-s}\alpha_-)$ is Liouville in $\mathbb{R}\times Y$. Consider also the \emph{Giroux $2\pi$-torsion domain} modeled on $(Y,\alpha_+,\alpha_-)$ given by the contact manifold $(GT,\xi_{GT}):=(Y\times [0,2\pi]\times S^1,\ker\lambda_{GT})$, where 
\begin{equation}
\begin{split}
\lambda_{GT}&=\frac{1+\cos(r)}{2}\alpha_++\frac{1-\cos(r)}{2}\alpha_-+\sin (r)d\theta\\
\end{split}
\end{equation} and the coordinates are $(r,\theta)\in [0,2\pi] \times S^1$. Say that a contact manifold $(M^{2n+1},\xi)$ has \emph{Giroux torsion} whenever it admits a contact embedding of $(GT,\xi_{GT})$. In this situation, denote by $\mathcal{O}(GT)\subseteq H^2(M;\mathbb{R})$ the annihilator of $\mathcal{R}_{GT}:=H_1(Y;\mathbb{R})\otimes H_1(S^1;\mathbb{R})$, viewed as a subspace of $H_2(M;\mathbb{R})$. The following was conjectured in \cite{MNW}, and proved in \cite{Mo}:

\begin{thm}\cite{Mo}\label{AlgvsGiroux}
If a contact manifold $(M^{2n+1},\xi)$ has Giroux torsion, then it has $\Omega$-twisted algebraic 1-torsion, for every $[\Omega] \in \mathcal{O}(GT)$, where $GT$ is a Giroux $2\pi$-torsion domain embedded in $M$. 
\end{thm}

Conversely, one could ask whether there exist examples of non-fillable contact manifolds without Giroux torsion. Examples of 3-dimensional weakly but not strongly fillable contact manifolds without Giroux torsion were given in \cite{NW11}. In higher dimensions, the following theorem can be proved without appealing to the abstract perturbation scheme for SFT (see Disclaimer \ref{disc} below). We use the fact that the unit cotangent bundle of a hyperbolic surface fits into a cylindrical semi-filling \cite{McD}.

\begin{thm}\label{nocob1} Let $(M_0^{5},\xi_0)$ be a $5$-dimensional contact manifold with Giroux torsion, and let $Y$ be the unit cotangent bundle of a hyperbolic surface. If $(M=\Sigma \times Y,\xi_k)$ is the corresponding $5$-dimensional contact manifold with $k\geq 3$, then there is no exact symplectic cobordism having $(M_0,\xi_0)$ as the convex end, and $(M,\xi_k)$ as the concave end.
\end{thm}

In particular, we obtain

\begin{cor}\label{corGT}
If $Y$ is the unit cotangent bundle of a hyperbolic surface, and $(M=\Sigma \times Y,\xi_k)$ is the corresponding $5$-dimensional contact manifold of Theorem \ref{thm5d} with $k\geq 3$, then $(M,\xi_k)$ does not have Giroux torsion.
\end{cor}

In the other direction of Theorem \ref{AlgvsGiroux}, examples of contact 3-manifolds which have $1$-torsion but not Giroux torsion where constructed in \cite{LW}. In higher dimensions, we have fairly strong geometric reasons for the following:

\begin{conj}\label{1tc} The examples of Corollary \ref{corGT} have untwisted algebraic 1-torsion (for \emph{any} $k\geq 1$).
\end{conj}

For $k\geq 3$, this would give an infinite family of contact manifolds with no Giroux torsion and algebraic 1-torsion in dimension $5$.

The interesting thing about this conjecture is that it builds on the relationship between SFT and string topology as discussed in \cite{CL09}. There are certain non-zero counts of punctured tori in the symplectization of $(Y,\alpha)$, where $\alpha$ is the standard Liouville form, which survive in $M$. They are given in terms of the coefficients of the Goldman--Turaev cobracket operation on strings in the underlying hyperbolic surface. One can find elements in the SFT algebra of $(Y,\alpha)$ whose differential in $\mathbb{R} \times Y$ has non-zero contributions from these tori and other pair of pants configurations, but by index considerations, only the former survive in $M$. This is how 1-torsion should arise. In fact, out of 35 possibilities, consisting of both cylinders and 1-punctured tori, this is the \emph{only} possible configuration from which 1-torsion can arise (see Section \ref{proooof}). We will therefore refer to it as \emph{sporadic}. This conjecture relies on technicalities regarding obstruction bundles. However, we provide a heuristic argument as to why we expect it to be true.

Theorem \ref{nocob1} would be then a result which is beyond the scope of algebraic torsion, since the presence of such a cobordism would only yield the already known fact that both ends have algebraic 1-torsion. In fact, in the proof of Theorem \ref{nocob1}, only holomorphic \emph{cylinders} play a role, whereas 1-punctured tori do not. In contrast, both are taken into consideration for $1$-torsion. This suggests the existence of an invariant more subtle than algebraic torsion, which we suspect would be obtained from the \emph{rational} SFT. Also, the computations in \cite{LW} could also be interpreted in this way.

Putting Theorem \ref{thm5d}, together with Corollary \ref{corGT}, we obtain the following:

\begin{cor}\label{corG}
There exist infinitely many non-diffeomorphic $5$-dimensional contact manifolds $(M,\xi)$ which are tight, not strongly fillable, and which do not have Giroux torsion.
\end{cor}

To our knowledge, there are no other known examples of higher-dimensional contact manifolds as in Corollary \ref{corG}. 

\vspace{0.5cm}

One can twist the contact structure of Theorem \ref{thm5d} close to the dividing set, by performing the \emph{$l$-fold Lutz--Mori twist} along a hypersurface $H$ lying in $\partial(\bigsqcup^k Y \times I \times S^1)$. This notion was defined in \cite{MNW}, and builds on ideas by Mori in dimension 5 \cite{Mori09}. It consists in gluing copies, along $H$, of a \emph{$2\pi l$-Giroux torsion domain} $(GT_l,\xi_{GT}):=(Y \times [0,2\pi l]\times S^1,\ker \lambda_{GT})$, the contact manifold obtained by gluing $l$ copies of $GT=GT_1$ together. The resulting contact structures are all homotopic as almost contact structures. By construction, all of these have Giroux torsion, so by Theorem \ref{AlgvsGiroux} they have $\Omega$-twisted $1$-torsion, for $[\Omega] \in \mathcal{O}=\mathcal{O}(GT)$.

As a corollary of Theorem \ref{nocob1}, we get:

\begin{cor}\label{nocob}
Let $Y$ be the unit cotangent bundle of a hyperbolic surface, and let $(M=\Sigma \times Y,\xi)$ be the corresponding $5$-dimensional contact manifold of Theorem \ref{thm5d}, with $k\geq 3$. If $(M,\xi_l)$ denotes the contact manifold obtained by an $l$-fold Lutz--Mori twist of $(M,\xi)$, then there is no exact symplectic cobordism having $(M,\xi_l)$ as the convex end, and $(M,\xi)$ as the concave end (even though the underlying manifolds are diffeomorphic, and the contact structures are homotopic as almost contact structures).  
\end{cor}

In the proof of Theorem \ref{nocob1}, we make use of obstruction bundles, in the sense of Hutchings--Taubes \cite{HT1,HT2}. We prove they exist in our setup, under the condition that they exist inside the leaves of a codimension-2 holomorphic foliation. As a byproduct, we derive a result for super-rigidity of holomorphic cylinders in $4$-dimensional symplectic cobordisms, which might be of independent interest. This is a natural adaptation of the results of Section 7 in \cite{Wen6} to the punctured setting. Recall that a somewhere injective holomorphic curve is super-rigid if it is immersed, has index zero, and the normal component of the linearized Cauchy--Riemann operator of every multiple cover is injective.  

\begin{thm}\label{superig} For generic $J$, every somewhere injective holomorphic cylinder in a $4$-dimensional symplectic cobordism, with index zero, and vanishing Conley--Zehnder indices (in some trivialization), is super-rigid.
\end{thm}

\begin{disclaimer}\label{disc} Since the statements of our results make use of machinery from Symplectic Field Theory, they come with the standard disclaimer that they assume that its analytic foundations are in place. However, we have taken special care in that the approach taken not only provides results that will be fully
rigorous after the polyfold machinery of Hofer--Wysocki--Zehnder is complete, but also gives several direct results that are \emph{already} rigorous.
\end{disclaimer}

\subparagraph*{Acknowledgements.} 

I thank my PhD supervisor, Chris Wendl, for introducing me to this project and for his support and patience throughout its duration. To Richard Siefring, for very helpful conversations and for co-authoring Appendix C in \cite{Mo2}. To Janko Latschev and Kai Cieliebak, for going through the long process of reading \cite{Mo2}. To Patrick Massot, Sam Lisi, Michael Hutchings and Momchil Konstantinov, for helpful conversations/correspondence on different topics. 

This research, forming part of the author's PhD thesis, has been partly carried out in (and funded by) University College London (UCL) in the UK, and by the Berlin Mathematical School (BMS) in Germany.

\section*{Guide to the document}

The main construction is carried out in Section \ref{RulingOut}, which corresponds to the building of our model contact manifolds in any odd dimension, from a ``dual'' point of view as in \cite{Mo}. We construct a foliation by holomorphic hypersurfaces in Section \ref{holhyps}. In Section \ref{curvesvshyp}, we use the results from \cite{MS} to show that relevant holomorphic curves lie in the leaves of the foliation. We discuss obstruction bundles in Section \ref{obbundles} from a fairly general point of view, and we prove their existence for curves in the symplectization of our models, assuming their existence leaf-wise. 

We illustrate our techniques in Section \ref{proooof}, where we study the algebraic $1$-torsion of a special family of contact manifolds in dimension $5$. We discuss the sporadic configuration in Section \ref{esp}, and its relationship to string topology. Theorem \ref{nocob1} is proved in Section \ref{5dmodel}. 
 
\vspace{0.5cm}
In Appendix \ref{obpants}, we derive Theorem \ref{superig}. In Appendix \ref{LMtwists} we describe the Lutz--Mori twists, and apply it to our model contact manifolds, to obtain the contact structures of Corollary \ref{nocob}. 

\newpage

\paragraph*{Basic notions}\label{basicn}

A contact form in a $(2n+1)$-dimensional manifold $M$ is a $1$-form $\alpha$ such that $\alpha \wedge d\alpha^{n}$ is a volume form, and the associated contact structure is $\xi=\ker \alpha$ (we will assume all our contact structures are co-oriented). The Reeb vector field associated to $\alpha$ is the unique vector field $R_\alpha$ on $M$ satisfying $$\alpha(R_\alpha)=1,\;i_{R_\alpha}d\alpha=0$$

A $T$-periodic Reeb orbit is $(\gamma,T)$ where $\gamma: \mathbb{R}\rightarrow M$ is such that $\dot{\gamma}(t)=TR_\alpha(\gamma(t))$, $\gamma(1)=\gamma(0)$. We will often just talk about a Reeb orbit $\gamma$ without mention to $T$, called its period, or action. If $\tau>0$ is the minimal number for which $\gamma(\tau)=\gamma(0)$, and $k \in \mathbb{Z}^+$ is such that $T=k\tau$, we say that the covering multiplicity of $(\gamma,T)$ is $k$. If $k=1$, then $\gamma$ is said to be simply covered (otherwise it is multiply covered). A periodic orbit $\gamma$ is said to be non-degenerate if the restriction of the time $T$ linearised Reeb flow $d\varphi^T$ to $\xi_{\gamma(0)}$ does not have $1$ as an eigenvalue. More generally, a Morse--Bott submanifold of $T$-periodic Reeb orbits is a closed submanifold $N \subseteq M$ invariant under $\varphi^T$ such that $\ker(d\varphi^T - \mathds{1})=TN$, and $\gamma$ is Morse--Bott whenever it lies in a Morse--Bott submanifold, and its minimal period agrees with the nearby orbits in the submanifold. The vector field $R_\alpha$ is non-degenerate/Morse--Bott if all of its closed orbits are non-degenerate/Morse--Bott.

A stable Hamiltonian structure (SHS) on $M$ is a pair $\mathcal{H}=(\Lambda,\Omega)$ consisting of a closed $2$-form $\Omega$ and a $1$-form $\Lambda$ such that 

$$\ker \Omega \subseteq \ker d\Lambda,\;\mbox{and}\;\Omega\vert_{\xi} \mbox{ is non-degenerate, where } \xi=\ker \Lambda$$ 

In particular, $(\alpha,d\alpha)$ is a SHS whenever $\alpha$ is a contact form. The Reeb vector field associated to $\mathcal{H}$ is the unique vector field on $M$ defined by $$\Lambda(R)=1,\; i_{R}\Omega=0$$ There are analogous notions of non-degeneracy/Morse--Bottness for SHS.

A symplectic form in a $2n$-dimensional manifold $W$ is a $2$-form $\omega$ which is closed and non-degenerate. A Liouville manifold (or an exact symplectic manifold) is a symplectic manifold with an exact symplectic form $\omega=d\lambda$, and the associated Liouville vector field $V$ is defined by the equation $i_{V}d\lambda=\lambda$. Any Liouville manifold is necessarily open. A boundary component $M$ of a Liouville manifold (endowed with the boundary orientation) is convex if the Liouville vector field is positively transverse to $M$, and is concave, if it is so negatively. An exact cobordism from a (co-oriented) contact manifold $(M_+,\xi_+)$ to $(M_-,\xi_-)$ is a compact Liouville manifold $(W,\omega=d\alpha)$ with boundary $\partial W=M_+\bigsqcup M_-$, where $M_+$ is convex, $M_-$ is concave, and $\ker \alpha\vert_{M_\pm}=\xi_\pm$. Therefore, the boundary orientation induced by $\omega$ agrees with the contact orientation on $M_+$, and differs on $M_-$. A Liouville filling (or a Liouville domain) of a --possibly disconnected-- contact manifold $(M,\xi)$ is a compact Liouville cobordism from $(M,\xi)$ to the empty set. A strong symplectic cobordism and a strong filling are defined in the same way, with the difference that $\omega$ is exact only in a neighbourhood of the boundary of $W$ (so that the Liouville vector field is defined in this neighbourhood, but not necessarily in its complement).

The symplectization of a contact manifold $(M,\xi=\ker \alpha)$ is the symplectic manifold $(\mathbb{R}\times M, \omega=d(e^a\alpha))$, where $a$ is the $\mathbb{R}$-coordinate. In particular, it is a non-compact Liouville manifold. Similarly, the symplectization of a stable Hamiltonian manifold $(M,\Lambda,\Omega)$ is the symplectic manifold $(\mathbb{R}\times M, \omega^\varphi)$, where $\omega^\varphi=d(\varphi(a)\Lambda)+\Omega$, and $\varphi$ is an element of the set $$\mathcal{P}=\{\varphi \in C^\infty(\mathbb{R},(-\epsilon,\epsilon)):\varphi^\prime>0\}$$ Here, $\epsilon>0$ is chosen small enough so that $\omega^\varphi$ is indeed symplectic. An $\mathcal{H}$-compatible (or simply cylindrical) almost complex structure on a symplectization $(W=\mathbb{R}\times M,\omega^\varphi)$ is $J \in \mbox{End}(TW)$ such that $$J \mbox{ is } \mathbb{R}\mbox{-invariant}, J^2=-\mathds{1},\;J(\partial_a)=R,\;J(\xi)=\xi,\;J\vert_{\xi} \mbox{ is }\Omega \mbox{-compatible }$$ The last condition means that $\Omega(\cdot,J\cdot)$ defines a $J$-invariant Riemannian metric on $\xi$. If $J$ is $\mathcal{H}$-compatible, then it is easy to check that it is $\omega^\varphi$-compatible, which means that $\omega^\varphi(\cdot,J\cdot)$ is a $J$-invariant Riemannian metric on $\mathbb{R}\times M$.

We will consider, for cylindrical $J$, punctured $J$-holomorphic curves $u:(\dot{\Sigma},j)\rightarrow (\mathbb{R}\times M,J)$ in the symplectization of a stable Hamiltonian manifold $M$, where $\dot{\Sigma}=\Sigma \backslash \Gamma$, $(\Sigma,j)$ is a compact connected Riemann surface, and $u$ satisfies the nonlinear Cauchy--Riemann equation $du \circ j = J \circ u$. We will also assume that $u$ is asymptotically cylindrical, which means the following. Partition the punctures into \emph{positive} and \emph{negative} subsets $\Gamma = \Gamma^+\cup \Gamma^-$, and at each $z \in \Gamma^\pm$, choose a biholomorphic identification of a punctured neighborhood of $z$ with the half-cylinder $Z_\pm$, where $Z_+ = [0, \infty) \times S^1$ and $Z_- = (-\infty, 0] \times S^1$. Then writing $u$ near the puncture in cylindrical coordinates $(s, t)$, for $|s|$
sufficiently large, it satisfies an asymptotic formula of the form 
$$u \circ\phi(s, t) = exp_{(Ts,\gamma(T t))} h(s, t)$$
Here $T > 0$ is a constant, $\gamma : \mathbb{R} \rightarrow M$ is a $T$-periodic Reeb orbit, the exponential map is defined with respect to any $\mathbb{R}$-invariant metric on $\mathbb{R} \times M$, $h(s, t) \in \xi_{\gamma(T t)}$ goes to $0$ uniformly in $t$ as $s \rightarrow \pm \infty$ and $\phi : Z_\pm \rightarrow Z_\pm$ is a smooth embedding such that $\phi(s, t) - (s + s_0, t + t_0) \rightarrow 0$ as $s \rightarrow \pm \infty$ for some constants $s_0 \in \mathbb{R}$, $t_0 \in S^1$. We will refer to punctured asymptotically cylindrical $J$-holomorphic curves simply as $J$-holomorphic curves.

Observe that, for any closed Reeb orbit $\gamma$ and cylindrical $J$, the trivial cylinder over $\gamma$, defined as $\mathbb{R}\times \gamma$, is $J$-holomorphic. 

The Fredholm index of a punctured holomorphic curve $u$ which is asymptotic to non-degenerate Reeb orbits in a $(2n+2)$-dimensional symplectization $W^{2n+2}=\mathbb{R} \times M$ is given by the formula 

\begin{equation} \label{index} 
ind(u)=(n-2)\chi(\dot{\Sigma})+2c_1^\tau(u^*TW)+\mu^\tau_{CZ}(u)
\end{equation} 

Here, $\dot{\Sigma}$ is the domain of $u$, $\tau$ denotes a choice of trivializations for each of the bundles $\gamma_z^*\xi$, where $z \in \Gamma$, at which $u$ approximates the Reeb orbit $\gamma_z$. The term $c_1^{\tau}(u^*TW)$ is the relative first Chern number of the bundle $u^*TW$. In the case $W$ is $2$-dimensional, this is defined as the algebraic count of zeroes of a generic section of $u^*TW$ which is asymptotically constant with respect to $\tau$. For higher-rank bundles, one determines $c_1^\tau$ by imposing that $c_1^\tau$ is invariant under bundle isomorphisms, and satisfies the Whitney sum formula (see e.g.\ \cite{Wen8}). The term $\mu^\tau_{CZ}(u)$ is the total Conley--Zehnder index of $u$, given by 

$$
\mu^\tau_{CZ}(u)=\sum_{z \in \Gamma^+}\mu^\tau_{CZ}(\gamma_z) - \sum_{z \in \Gamma^-}\mu^\tau_{CZ}(\gamma_z)
$$ 

\section{Model contact manifolds: a dual point of view.} 
\label{RulingOut}

We will construct a contact model for the underlying manifold $M=Y \times \Sigma$. Here, as described in the introduction, $Y$ fits into a cylindrical semi-filling, and $\Sigma$ is an orientable genus $g$ surface, decomposed into a genus zero piece, and a positive genus piece along a dividing set of $k$ circles. The main features of this contact model are: closed Reeb orbits of low action correspond to pairs $(p,\gamma)$ of critical points $p$ of suitable Morse functions on $\Sigma_\pm$, and closed $R_\pm$-orbits $\gamma$ in $(Y,\alpha_\pm)$; and we will have a foliation of its symplectization by holomorphic hypersurfaces, for a suitable almost complex structure compatible with a SHS deforming the contact data. These project to flow lines on the surface $\Sigma$, and they can be identified with either cylindrical completions of the Liouville domain $Y\times I$ or symplectizations of its boundary components. The results in \cite{MS} may then be used to prove that any holomorphic curves necessarily lies in a hypersurface of this foliation, which is the key tool for restricting their behaviour. After constructing the model, we will restrict our attention to a specific family of $5$-dimensional models for which we prove Theorem \ref{nocob1}.  
 
\subsection{Construction of the model contact manifolds} 

Consider $Y^{2n-1}$ such that $(Y\times I,d\alpha)$ is a cylindrical Liouville semi-filling. We assume that $\alpha=\{\alpha_r\}_{r\in I}$ is a 1-parameter family of 1-forms such that $\alpha_r$ is contact for $r\neq 0$. All known examples of such semi-fillings satisfy this condition. We will denote by $\xi_\pm=\ker \alpha_\pm$ the contact structures on $Y_\pm:=Y\times \{\pm 1\}$, where $\alpha_\pm=\alpha_{\pm 1}=\alpha\vert_{Y_\pm}$, by $V$ the Liouville vector field associated to $d\alpha$, by $R_r$ the Reeb vector field associated to $\alpha_r$ for $r\neq 0$, and $R_\pm:=R_{\pm 1}$. Consider also $\Sigma=\Sigma_-\cup_\partial \Sigma_+$ a genus $g$ surface obtained by gluing a genus zero surface $\Sigma_-$ to a genus $g-k+1>0$ surface $\Sigma_+$, along $k$ boundary components in an orientation-reversing manner, and the manifold $M=Y\times \Sigma$. Throughout the paper, we shall make the convention that whenever we deal with equations involving $\pm$'s and $\mp$'s, one has to interpret them as to having a different sign according to the region (the ``upper'' sign denotes the ``plus'' region, and the ``lower'', the ``minus'' region).

Take collar neighbourhoods of each boundary component of $\Sigma_\pm$ of the form $\mathcal{N}(\partial \Sigma_\pm):=[-1,0]\times S^1 $, with coordinates $(t_\pm,\theta) \in [-1,0] \times S^1$, such that $\partial\Sigma_\pm=\{t_\pm=0\}$. Glue these surfaces together along a cylindrical region $\mathcal{U}$, having $k$ components, each of which looks like $[-1,1]\times S^1 \ni (r,\theta)$, so that $r=\mp t_\pm$. We now use the fact that each surface $\Sigma_\pm$ carries a \emph{Stein} structure providing a Stein filling of its 1-dimensional convex contact boundary. That is, we may take Morse functions $h_\pm$ on $\Sigma_\pm$ which are plurisubharmonic with respect to a suitable complex structure $j_\pm$, which in the collar neighbourhoods look like $h_\pm(t_\pm,\theta)=e^{t_\pm}$ and $j_\pm(\partial_{t_\pm})=\partial_\theta$. Moreover, we have that $\lambda_\pm=d^{\mathbb{C}}h_\pm=-dh_\pm \circ j_\pm$ is a Liouville form, with Liouville vector field given by $\nabla h_\pm$, where the gradient is computed with respect to the metric $g_\pm=g_{d\lambda_\pm,j_\pm}=d\lambda_\pm(\cdot,j_\pm \cdot)$ (see Lemma 4.1 in \cite{LW} for details). Thus, in the collar neighbourhoods, we have $\lambda_\pm = e^{t_\pm}d\theta=e^{\mp r}d\theta$ and $\nabla h_\pm = \partial_{t_\pm}=\mp \partial_r$. The Hamiltonian vector field $X_{h_\pm}=j_\pm \nabla h_\pm$, computed with respect to the symplectic form $d\lambda_\pm$ inducing the orientation in $\Sigma_\pm$, is tangent to the contact-type level sets. 

We take the orientation on $M= Y\times \Sigma$ to be the one induced by the $\alpha_-$ and $d\lambda_-$ on each respective factor. We take both $h_\pm$ to have a unique minimum in $\Sigma_\pm$, and we will make the convention of calling the minimum of $h_+$, the \emph{maximum} (cf. Remark \ref{admissible} below).

We define, for $K_+>K_-\gg 0$ large constants, a function $g_0:[-1,1] \rightarrow \mathbb{R}$, such that it is constant equal to $K_-$ in $[-1,-1+\delta)$ for some small $\delta>0$, equal to $K_+$ in $(1-\delta,1]$, $g_0> 0$ everywhere, and $g_0^\prime\geq 0$ with strict inequality in the interval $(-1+\delta,1-\delta)$. The parameter $\delta>0$ is chosen so that the Liouville vector field in $Y \times [-1,1]$ is given by $V=\pm \partial_r$ in the corresponding components of $\{|r|>1-\delta\}$.  Define also a function $\gamma: [-1,1]\rightarrow \mathbb{R}$ so that $\gamma(r)=- e^{\mp r}$ in $\{|r|>1-\delta\}$, and $\gamma<0$. Take $$g_\epsilon(r)=g_0(r)+\epsilon^2\gamma(r),$$ where $\epsilon$ is chosen small enough so that $g_\epsilon>0$.

Now define 
$$h_\epsilon=\left\{\begin{array}{ll} K_\pm - \epsilon^2 h_\pm,& \mbox{ in } \Sigma\backslash \mathcal{U} \\
g_\epsilon,& \mbox { in } \mathcal{U} 
\end{array}\right.
$$

This is a (well-defined and smooth) function on $\Sigma$, and $\epsilon$ can be chosen small enough so that $h_\epsilon>0$.

\vspace{0.5cm}

Next, choose a smooth function $\psi:[-1,1]\rightarrow \mathbb{R}^+$ satisfying $\psi(r)=e^{\mp r}$ in $\{|r|>1-\delta/3\}$, $\psi\equiv 1$ on $[-1+\delta,1-\delta]$, $r\psi^\prime(r)\leq 0$. Define $$f_\epsilon(r)=\epsilon\psi(r),$$ so that $f_\epsilon>0$, $f_\epsilon^\prime\equiv 0$ along $[-1+\delta,1-\delta]$. 

With these choices, let $\Lambda_\epsilon$ be the 1-parameter family of 1-forms in $M=Y \times \Sigma$ given by 
\begin{equation}
\label{modelBform}
\Lambda_\epsilon=\left\{ \begin{array}{ll}\epsilon\lambda_\pm + h_\epsilon \alpha_\pm, &\mbox{ in } Y\times (\Sigma\backslash\mathcal{U})\\
f_\epsilon d\theta + h_\epsilon \alpha, &\mbox{ in } Y\times \mathcal{U}
\end{array}\right.
\end{equation}

We can slightly modify $\alpha$ so that the different expressions in $\Lambda_\epsilon$ glue together smoothly: we have $\alpha=e^{\pm r-1}\alpha_\pm$ in the corresponding components of the region $\{|r|\geq 1-\delta\}$), and we can replace $\alpha$ by a 1-form (of the same name) which looks like $\phi(r)\alpha_\pm$. Here, $\phi: \{|r|\geq 1-\delta\}\rightarrow \mathbb{R}$ is a smooth function which coincides with $e^{\pm r-1}$ near $|r|=1-\delta$, equals $1$ on $|r|=1$, and its derivative is non-negative/non-positive in the positive/negative components of $\{|r|\geq 1-\delta\}$.   

We shall refer to this family of contact forms as \emph{model B}, in contrast with \emph{model A}, which is the contact form for $M$ constructed in \cite{Mo}.

\begin{lemma} \label{contact} $\;$
\begin{enumerate}
\item For a fixed $K_+\gg 0$, $\Lambda_\epsilon$ is a contact form for $K_-$ sufficiently close to $K_+$, and sufficiently small $\epsilon>0$.
\item Model A is isotopic to model B.

\end{enumerate}
\end{lemma}

\begin{proof} The first claim is straightforward to check. For the second, see \cite{Mo2}.
\end{proof}

\subsection{Deformation to a SHS along the cylindrical region}

To obtain a SHS deforming the contact data, which we will need to construct holomorphic hypersurfaces, we homotope the Liouville vector field over $\{\rho \in [-1,1]\}$, as follows. Choose a bump function $\kappa:\mathbb{R}\rightarrow \mathbb{R}$ which equals $1$ outside of the unit interval, equals zero in $[-1+\delta,1-\delta]$, and sign$(\kappa^\prime(\rho))=$sign$(\rho)\neq 0$ for $\rho \in \{1-\delta<|\rho|<1\}$.

Define
$$
W^s_\epsilon:=\left\{\begin{array}{ll}Z_\epsilon, &\mbox{ in }E^{\infty,\infty}\backslash\left( \{|\rho|\leq 1\}\cap E^{\infty,\infty}(t)\right)\\
X +(s\kappa(\rho)+1-s)V_\epsilon, &\mbox{ in } \{|\rho|\leq 1\}\cap E^{\infty,\infty}(t) 
\end{array}\right.
$$

One can check that 
$$
\mathcal{H}_\epsilon^s:=(\Lambda_\epsilon^s,\Omega_\epsilon)=(i_{W_s^\epsilon}d\Lambda_\epsilon\vert_{TM}, d\Lambda_\epsilon\vert_{TM})_{s\in[0,1]}
$$
is a family of SHS's on $M$, which deform model B, and such that $\Lambda_\epsilon^s$ is contact for $s<1$. 

\vspace{0.5cm}

Focusing now in the case $L=1$, $Q=0$, we have $\Lambda_\epsilon^s=\epsilon \lambda+(s\kappa+1-s)h_\epsilon\alpha$, which we can explicit write as
$$
\Lambda_\epsilon^s=\left\{\begin{array}{ll} \epsilon \lambda_\pm + h_\epsilon\alpha_\pm,& \mbox{ in } M_P^\pm=Y\times (\Sigma\backslash \mathcal{U})\\
f_\epsilon(r)d\theta + (s\kappa(F_\pm(r)\pm 1)+1-s)h_\epsilon\alpha,& \mbox{ in } M_C^\pm=Y\times\left(\mathbb{R}^\pm\cap \{|r|>1-\delta\}\right)\times S^1\\
\epsilon d\theta+(1-s)h_\epsilon \alpha, &\mbox{ in } M_Y= Y \times [-1+\delta,1-\delta]\times S^1,
\end{array}\right.
$$
where we denote $f_\epsilon(r)=\epsilon e^{G_\pm(r)}$ (cf.\ (\ref{modelBform})).

\vspace{0.5cm}

One may check that the Reeb vector field $R_\epsilon^s$ associated to $\mathcal{H}^s_\epsilon$ is given by
$$
R_\epsilon^s=\left\{\begin{array}{ll} R_\epsilon,& \mbox{ in } M_P^\pm\\ 
\frac{1}{D_\epsilon^s +f_\epsilon h_\epsilon \frac{\partial \alpha_r}{\partial r}(R_r)}((h_\epsilon^\prime + h_\epsilon \frac{\partial \alpha_r}{\partial r}(R_r))\partial_\theta-f_\epsilon^\prime R_r),& \mbox{ in } M_C^\pm\\
\partial_\theta/\epsilon, &\mbox{ in } M_Y
\end{array}\right. 
$$ 
where $D_\epsilon^s:= f_\epsilon h_\epsilon^\prime - (s\kappa(F_\pm \pm 1)+1-s)h_\epsilon f_\epsilon^\prime$. 

Every pair $(p,\gamma)$, consisting of a critical point $p$ of $h_\pm$ and a closed $R_\pm$-orbit $\gamma$, gives rise to a closed Reeb orbit of $R_\epsilon^s$, of the form $\gamma_p:=(\gamma,p) \in Y \times \Sigma$. The closed Reeb orbits which do not correspond to critical points of $h_\pm$ can be made to have arbitrarily large period by taking $\epsilon$ sufficiently close to zero. Therefore we can find an action threshold $T_\epsilon>0$ ($s$-indfependent), such that $\lim_{\epsilon \rightarrow 0} T_\epsilon = +\infty$, and such that every closed Reeb orbit of $R_\epsilon^s$ of action less than $T_\epsilon$ either lies in $Y \times \mathcal{U}$, or is a cover $\gamma_p^l$ of a simple closed Reeb orbit $\gamma_p$, for $l\leq N_\epsilon$, where $N_\epsilon:=\max\{n\in\mathbb{N}:\gamma^l \mbox{ has action less than } T_\epsilon\}$. 

The following lemma will be used in the proof of Theorem \ref{nocob1}.

\begin{lemma}\label{sympcob}
The 2-form $\omega=d(e^\epsilon \Lambda_\epsilon^s)$ is symplectic in $(0,\epsilon_0]\times M$, for sufficiently small $\epsilon_0>0$ and for $0\leq s< 1$, and induces the orientation given by the natural product orientation (where the orientation on $M$ is as described before Lemma \ref{contact}).
\end{lemma}
\begin{proof} It follows by direct computations \cite{Mo2}.
\end{proof}

\subsection{Compatible almost complex structure} 
\label{compJmodelB}

We proceed now to define a $\mathcal{H}_\epsilon^s$-compatible (and non-generic) almost complex structure $J_\epsilon^s$ in $\mathbb{R}\times M$. We define it on $\xi_\epsilon^s=\ker \Lambda_\epsilon^s$ and extend it in a cylindrical way. 

We write 
\begin{equation}\label{xiB}
\xi^s_\epsilon=\left\{\begin{array}{ll} \xi_\pm \oplus L^s_\epsilon, &\mbox{ in } M_P^\pm \cup M_C^\pm\\
\widetilde{T(Y\times I)},& \mbox{ in } M_Y
\end{array}\right.
\end{equation}
where $\widetilde{T(Y\times I)}=\{v-(1-s)h_\epsilon\frac{\alpha(v)}{\epsilon}\partial_\theta:v \in T(Y\times I)\}$, and $L^s_\epsilon$ is a bundle of (real) rank 2. The latter may be written as 
\begin{equation}
L^s_\epsilon = \left\{\begin{array}{ll}
T\Sigma_\pm,& \mbox{ in crit}(h_\pm)\subseteq M_P^\pm\\ 
\left\langle \nabla h_\pm, h_\epsilon X_{h_\pm}-\epsilon |\nabla h_\pm|^2 R_\pm \right \rangle,& \mbox{ in } M_P^\pm\backslash \mbox{crit}(h_\pm)\\
\langle V=\pm \partial_r, W \rangle & \mbox{ in } M_C^\pm
\end{array}\right.
\end{equation}
where 
$$
W:=\mp f_\epsilon R_r\pm(s\kappa(F_\pm\pm 1)+1-s)h_\epsilon \partial_\theta
$$

Over the region $M_P^\pm=Y\times (\Sigma_\pm\backslash\mathcal{U})$, the differential of the projection $\pi: Y \times (\Sigma_\pm\backslash\mathcal{U})\rightarrow \Sigma_\pm\backslash\mathcal{U}$ is an isomorphism for every $\epsilon$ when restricted to $L^s_\epsilon$. We may then define  
\begin{equation}\label{splitJ}
J^s_\epsilon\vert_{\xi^s_\epsilon}= J_\pm \oplus (\pi\vert_{L^s_\epsilon})^*j_\pm,
\end{equation} with respect to the splitting above, where $J_\pm$ is a $d\alpha_\pm$-compatible almost complex structure in $\xi_\pm$, and $j_\pm$ is a $d\lambda_\pm$-compatible almost complex structure on $\Sigma_\pm$ which satisfies $j_\pm(\partial_{t_\pm})=\partial_\theta$ on the collar neighbourhoods.   

\vspace{0.5cm} 

Choose now a $d\alpha$-compatible almost complex structure $J_0$ in the Liouville domain  $Y\times I$ which is cylindrical in the cylindrical ends $\{|r|>1-\delta\}$. We impose that, along these, its restriction to $\xi_r$ coincides with $J_\pm$ . 

\vspace{0.5cm} 

Over the region $Y \times \mathcal{U}=Y \times I \times S^1,$ the projection $\pi_Y: Y \times \mathcal{U}\rightarrow Y \times I$ gives an isomorphism $$d\pi_Y\vert_{\xi^s_\epsilon}:\xi_\epsilon^s=\widetilde{T(Y \times I)} \stackrel{\simeq}{\longrightarrow} T(Y \times I),$$ and thus $J^S_\epsilon=\pi_Y^* J_0$ is an almost complex structure on $\xi^s_\epsilon$. 

\vspace{0.5cm}

In order to glue the two definitions along the region $\mathbb{R}\times M_C^\pm$, one computes that
$J^s_\epsilon(W)=s\kappa(F_\pm \pm 1)+1-s)h_\epsilon \partial_r,$ close to $r=\pm1$, and $J^s_\epsilon(W)=f_\epsilon \partial_r$ close to $r=\pm 1\mp \delta$. One may then define $J^s_\epsilon(W)=\beta^s_{\epsilon,\pm}(r) \partial_r$ for suitable interpolating functions $\beta^s_{\epsilon,\pm}(r)$. 

We thus get a well-defined almost complex structure $J^s_\epsilon$ over the whole model.

\paragraph*{Compatibility} One can check that $J_\epsilon^s$ is $\mathcal{H}_\epsilon^s$-compatible by explicit computation on the above splittings.

\begin{remark}\label{c1rk}
Observe that the splitting (\ref{xiB}) is holomorphic. One can directly check that
$c_1(\xi_\epsilon^1)=0$, whenever $c_1(\xi_\pm)=0$ (which holds for all known examples of cylindrical semi-fillings).
\end{remark}
  
\subsection{Foliation by holomorphic hypersurfaces} \label{holhyps}

The goal of this section is to construct a foliation of $\mathbb{R}\times M$ by holomorphic hypersurfaces (where the almost complex structure is the one compatible with the \emph{non-contact} SHS $\mathcal{H}_\epsilon^1$). 

\begin{prop}\label{curves} For $\epsilon \in (0,\epsilon_0]$, there exists a foliation $\mathcal{F}=\bigcup_{a,\nu} H_a^\nu$ of the symplectization of $\mathcal{H}_\epsilon^1$ by $J^1_\epsilon$-holomorphic hypersurfaces, which come in three types (see picture (\ref{fol2})):
\begin{enumerate}
\item (\textbf{cylindrical} hypersurfaces over critical points) For $p \in \mbox{crit}(h_\pm)$, there exists a \emph{cylindrical} hypersurface of the form $\mathbb{R}\times H_p$, where $H_p:=Y \times \{p\}$, a copy of the symplectization of $(Y,\alpha_\pm)$.

\item (positive/negative \textbf{flow-line} hypersurfaces contained in one side of the dividing set) These are parametrized by $$H_a^\nu=\{(a(s),\nu(s)): s \in \mathbb{R}\} \times Y,$$ where $a:\mathbb {R}\rightarrow \mathbb{R}$ is a proper function, and $\nu:\mathbb{R}\rightarrow \mbox{int}(\Sigma_\pm)$ is a non-constant negative/positive reparametrization of a flow-line for $h_\pm$.

If $\lim_{s \rightarrow \pm\infty}\nu(s)=p_\pm$, then $H_a^\nu$ is asymptotically cylindrical over $H_{p_\pm}$ (see \cite{Mo2, MS} for a definition), and has exactly one positive and one negative end. If the side containing it is the positive one (and thus $\nu$ connects the maximum to a hyperbolic point), the positive end corresponds to the maximum. If it is $\Sigma_-$ (and thus $\nu$ connects the minimum to a hyperbolic point), it corresponds to the minimum.   

\item (\textbf{page-like} hypersurfaces crossing sides of the dividing set) These may be written as
$$H_a^\nu=H_a^{\nu,-}\bigcup H_a^{\nu,cst}\bigcup H_a^{\nu,+},$$
where 
$$
H_a^{\nu,cst}=\{a\}\times Y \times [-1,1] \subseteq \mathbb{R}\times M_Y,
$$
$$
H_a^{\nu,+}:=\{(a_+(s),\nu_+(s)):s\in[0,+\infty)\}\times Y\subseteq \mathbb{R}\times M_P^+
$$
$$
H_a^{\nu,-}:=\{(a_-(s),\nu_-(s)):s\in(-\infty,0]\}\times Y\subseteq \mathbb{R}\times M_P^-
$$ Here, $a_+:[0,+\infty)\rightarrow \mathbb{R}$ and $a_-:(-\infty,0]\rightarrow \mathbb{R}$ are proper functions, and $\nu= \nu_+ \sqcup\nu_-$, for $\nu_+: [0,+\infty)\rightarrow \Sigma_+$ a non-constant reparametrization of a negative Morse flow line of $h_+$, and $\nu_-: (-\infty,0]\rightarrow \Sigma_-$ a non-constant reparametrization of a positive Morse flow line of $h_-$, both satisfying $\nu_\pm(0) \in \partial \Sigma_\pm$. If $\lim_{s \rightarrow \pm \infty}\nu_\pm(s)=p_\pm$, the hypersurface $H_a^\nu$ is asymptotically cylindrical over $H_{p_\pm}$, having exactly two positive ends.    
\end{enumerate}

\begin{figure}[t]\centering
\includegraphics[width=0.65\linewidth]{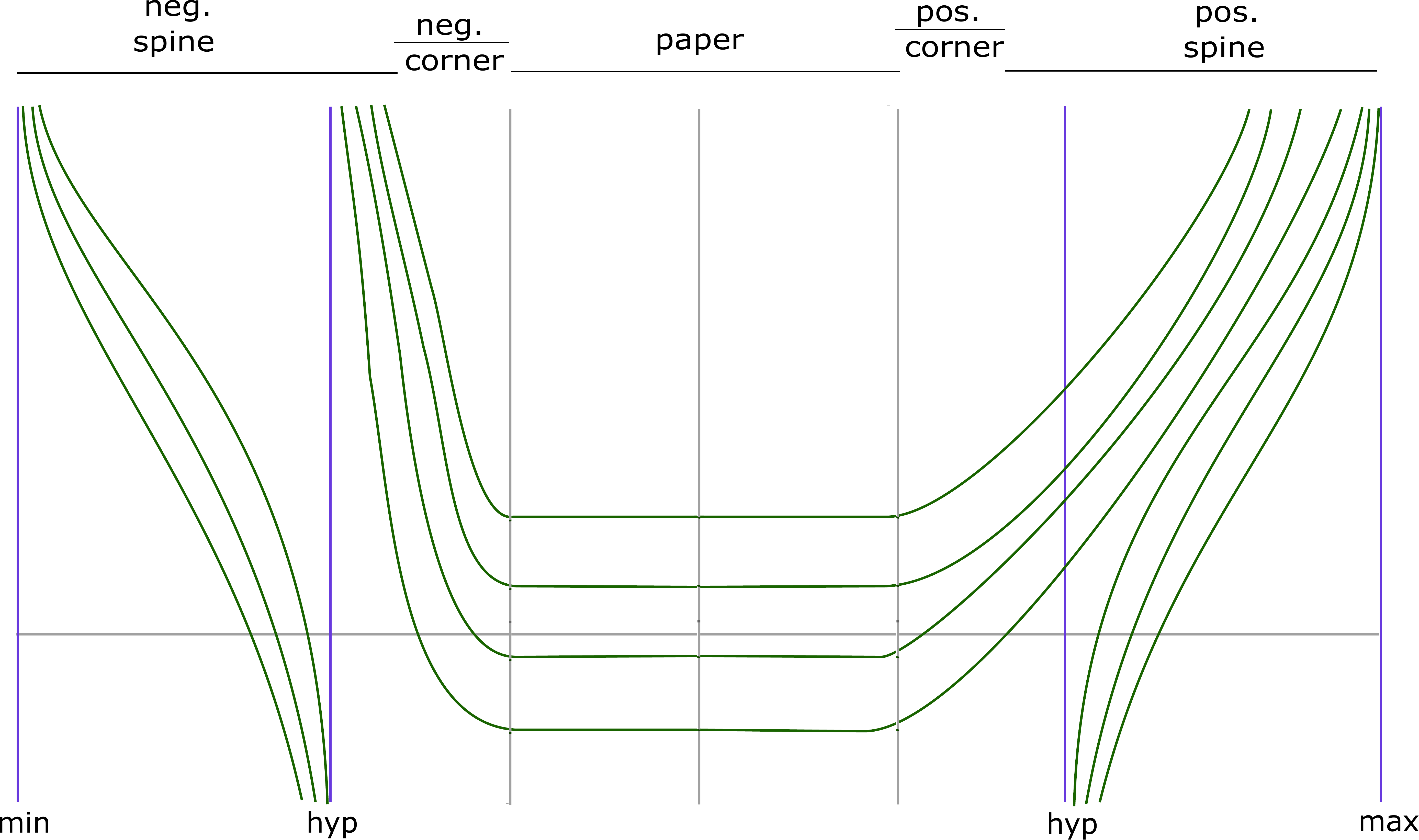}
\caption{\label{fol2} Foliation by holomorphic hypersurfaces. The vertical lines are the cylindrical hypersurfaces corresponding to critical points $min$ (minimum), $hyp$ (hyperbolic) or $max$ (maximum), and the non-vertical ones are the non-cylindrical $H_a^\nu$'s.}
\end{figure}
\end{prop}

\begin{proof}
One computes, along the paper, that  
$$
J_\epsilon^1 R_\pm=F\nabla^{g_\pm} h_\pm + G\partial_a,
$$
where
$$
F=\frac{\epsilon h_\epsilon}{h_\epsilon +\epsilon^2|\nabla h_\pm|^2},\; G=-h_\epsilon
$$
Observe that $F>0$, $G<0$.
Then
$$\langle R_\pm, J^1_\epsilon R_\pm\rangle \oplus \xi_\pm =\left\langle F\nabla^{g_\pm} h_\pm+G\partial_a\right\rangle \oplus TY$$ a $J^1_\epsilon$-invariant and integrable distribution over $M_P^\pm$. Integrating this distribution, one gets the hypersurfaces $H_a^\nu$ as in the statement. For instance, for the hypersurfaces of the second type, we take $a$ and $\nu$ to satisfy
$$\dot{\nu}(s)=\mp F(\nu(s))\nabla^{g_\pm} h_\pm(\nu(s))$$
$$\dot{a}(s)=\mp G(\nu(s))$$ 
For hypersurfaces of the third type, one can check that for $a_\pm$ and $\nu_\pm$ satisfying the same ODEs we can glue each piece as in the statement, to get a hypersurface $H_a^\nu$ whose tangent space is 
$$
TH_a^\nu=\left\{\begin{array}{ll} \langle R_\pm,J_\epsilon^1 R_\pm\rangle \oplus \xi_\pm,& \mbox{ over } M_P^\pm\cup  M_C^\pm\\
T(Y \times I), &\mbox{ over } M_Y
\end{array}\right.,
$$
which is clearly $J^1_\epsilon$-invariant. The assertion about the asymptotic behaviour of the ends of $H_a^\nu$ is deduced by the fact that orientations change as we move from one side to the other. The condition that the holomorphic hypersurfaces $H_a^\nu$ are asymptotically cylindrical follows from the fact that they project to $\Sigma$ as the Morse flow line $\nu$ \cite{Mo2}.    
\end{proof}

\begin{remark}\label{symhyp}
For a flow-line hypersurface $H_a^\nu$ which does not cross the dividing set, and $p_\pm=\lim_{s\rightarrow \pm\infty}\nu(s)$ a critical point, the projection map  
$$P_a^\nu: H_a^\nu=\{(a(s),\nu(s)):s\in \mathbb{R}\}\times Y \rightarrow \mathbb{R}\times H_{p_\pm}$$
$$(a(s),\nu(s),y)\mapsto  (s,y,p_\pm)$$ 
gives an identification of $H_a^\nu$ with $\mathbb{R}\times H_{p_\pm}$, which, up to replacing $s$ by $-s$ over $M_P^-$, is a biholomorphism. In particular it maps holomorphic curves in $H_a^\nu$ to holomorphic curves in $\mathbb{R}\times H_{p_\pm}$. The hypersurfaces which do cross sides are cylindrical completions of the Liouville domain $Y \times I$, and their cylindrical ends are identified with the symplectization of $(Y,\alpha_\pm)$ as above.  
\end{remark}

\begin{remark}\label{Weinstein} Assume that $(Y,\alpha_\pm)$ both satisfy the Weinstein conjecture, so that there are closed $\alpha_\pm$-Reeb orbits. Then, the hypersurfaces $H_a^\nu$ which stay on the same side of the dividing set contain finite energy $J^s_\epsilon$-holomorphic curves (for \emph{every} $s$). These are of the form
$$u^a_\gamma: \mathbb{R} \times S^1 \rightarrow \mathbb{R} \times M$$
$$u^a_\gamma(s,t)= (a(s),\nu(s),\eta(t)),$$ where $\nu$ and $a$ are as above, and $\eta$ is a \emph{closed} $\alpha_\pm$-Reeb orbit, satisfying $\lim_{s\rightarrow \pm\infty}\nu(s)=p_\pm$ for $p_\pm$ critical point. The asymptotic behaviour of $u_\gamma$ is thus the same as the hypersurface $H_a^\nu$ containing it. These cylinders correspond to trivial cylinders under the identification of the previous Remark \ref{symhyp}.

In general, there is no reason why holomorphic cylinders crossing sides should exist, and, in fact, there are examples for which they don't (see Section \ref{5dmodel}). 

\end{remark}

\begin{remark}\label{admissible}
In the language and notation of \cite{Mo2,MS}, the holomorphic foliation $\mathcal{F}$ is \emph{compatible} with the $H_\epsilon$-admissible fibration $(\Sigma,\pi_\Sigma,Y,H_\epsilon)$. Here, $\pi_\Sigma:M=Y \times \Sigma \rightarrow \Sigma$ is the natural projection, and $H_\epsilon$ is a suitably defined Morse function on $\Sigma$ having the same critical points as $h_\pm$, but up to orientation reversal on $\Sigma_+$, so that its unique maximum is the minimum of $h_+$ (see \cite{Mo2}). The \emph{binding} of $\mathcal{F}$ is $\bigcup_{p \in \crit(H_\epsilon)}H_p$, where each $H_p$ is a \emph{strong stable hypersurface}. The \emph{sign function} $\sign:\crit(H_\epsilon)\to\br{-1, 1}$ is defined by $\sign(p)=\pm 1$, if $p \in \crit(h_\pm)\subseteq \Sigma_\pm$.

We will use this fact in Section \ref{curvesvshyp}, together with the results from \cite{Mo2,MS}, to show that every holomorphic curve in $\mathbb{R}\times M$ with positive asymptotics corresponding to critical points has to lie in a leaf of $\mathcal{F}$.

\end{remark}

\subsection{Index Computations}

For closed Reeb orbits $\gamma_p: S^1 \rightarrow M$ of the form $\gamma_p(t)=(\gamma(t),p)$, for $p \in $ crit$(h_\pm)$ and $\gamma$ a closed $\alpha_\pm$-Reeb orbit, we have a natural identification 
\begin{equation}\label{ident}
\Gamma(\gamma_p^*\xi^s_\epsilon)=\Gamma(\gamma^*\xi_\pm)\oplus \Gamma(L^s_\epsilon\vert_p) = \Gamma(\gamma^*\xi_\pm)\oplus \Gamma(T_p\Sigma_\pm),
\end{equation}
where the second summand denotes the space of sections of the trivial line bundle over $S^1$ with fiber $L^s_\epsilon\vert_p=T_p\Sigma_\pm$. Therefore any trivialization $\tau$ of $\gamma^*\xi_\pm$ extends naturally to a trivialization of $\gamma_p^*\xi^s_\epsilon$, which we shall also denote by $\tau$. 

\begin{prop}\label{indexthm} 

Consider a trivialization $\tau$ of $\gamma^*\xi_\pm$ over the simply covered $\alpha_\pm$-Reeb orbit $\gamma$, inducing a natural trivialization $\tau^l$ of $(\gamma^l)^*\xi_\pm$ for every $l$. Then for each sufficiently small $\epsilon>0$ there exists a covering threshold $N_\epsilon$, satisfying $\lim_{\epsilon \rightarrow 0}N_\epsilon=+\infty$, such that the Conley--Zehnder index of $\gamma^l_p$ for $l\leq N_\epsilon$ with respect to the induced trivialization of $(\gamma_p^l)^*\xi^s_\epsilon$ given by the identification (\ref{ident}) is 

$$
\mu_{CZ}^{\tau^l}(\gamma^l_p)=\left\{
\begin{array}{ll}
\mu_{CZ}^{\tau^l}(\gamma^l)+1,& \mbox{ if } p \mbox{ is the maximum or minimum}\\
\mu_{CZ}^{\tau^l}(\gamma^l),& \mbox{ if } p \mbox{ is hyperbolic}
\end{array} \right.
$$

Here, $\mu_{CZ}^\tau(\gamma)$ denotes the Conley--Zehnder index of $\gamma$ in the case this Reeb orbit is non-degenerate or Morse--Bott.

\begin{remark}
In the Morse--Bott case, one needs to add suitable weights adapted to the spectral gap of the corresponding operator to obtain $\mu_{CZ}^\tau(\gamma)$, as explained e.g.\ in \cite{Wen1}.
\end{remark}

\end{prop}

\begin{proof} See \cite{Mo2}. 
\end{proof}

\subsection{Holomorphic curves lie in hypersurfaces}
\label{curvesvshyp}

In this section, we shall make use of intersection theory for punctured holomorphic curves and holomorphic hypersurfaces, as outlined in \cite{Mo2, MS}, to which we refer for the relevant definitions and notation. 

\vspace{0.5cm}

If we assume by perturbing that $(Y,\alpha_\pm)$ are non-degenerate contact manifolds, then we are in the situation of \cite{Mo2, MS}, as can be easily checked. For instance, one can compute that the splitting of the asymptotic operator $\mathbf{A}_{\gamma_p^l}$ into tangent and normal components is given by $$\mathbf{A}_{\gamma^l_p}=\mathbf{A}_{\gamma^l} \oplus \mathbf{A}_{\epsilon,\gamma_p^l}^\pm,$$
where the normal operator acting on $\Gamma(T_p\Sigma)$ is  
$$\mathbf{A}_{\epsilon,\gamma_p^l}^\pm=-J_0 \nabla_t - \epsilon T_{\gamma^l}\nabla^2_p h_\pm,$$ which has Conley--Zehnder index $\mu_{CZ}(\mathbf{A}_{\epsilon,\gamma_p^l}^\pm)=1-\mbox{ind}(h_\pm)=\sign(p)(\mbox{ind}_p(H_\epsilon)-1),$ where $H_\epsilon$ and the sign function are the ones of Remark \ref{admissible}.

\begin{prop}
\label{uniqhyp}
Suppose $u:\dot{S}\rightarrow \mathbb{R}\times M$ is a finite energy $J^1_\epsilon$-holomorphic curve which has all of its positive ends asymptotic to Reeb orbits of the form $\gamma_p^l$, with $l\leq N_\epsilon$. Then the image of $u$ is contained in a  leaf of the foliation $\mathcal{F}$.
\end{prop}
\begin{proof}
Given such a curve $u$, by Stokes' theorem we get that its negative ends have action bounded by $T_\epsilon$, and so also correspond to critical points. Since all of its asymptotics project to $\Sigma$ as points, we may define a map $v:=\pi_\Sigma \circ u: S\rightarrow \Sigma$ between \emph{closed} surfaces. We are therefore in the situation of \cite{MS}, which, since the map $\sign$ is surjective, implies the result.
\end{proof}

\begin{remark} Proposition \ref{uniqhyp} reduces the study of $J_\epsilon^1$-holomorphic curves/buildings inside $\mathbb{R}\times M$ to the study of $J_0$-curves/buildings inside the completion $W_0=\widehat{Y\times I}$ of the cylindrical Liouville semi-filling $Y\times I$. Recall that $J_0$ is any cylindrical $d\alpha$-compatible almost complex structure in $Y\times I$, and $W_0$ is obtained by attaching cylindrical ends to $\partial(Y\times I)=(Y_-,\alpha_-)\bigsqcup(Y_+,\alpha_+)$. Holomorphic buildings inside $W_0$ are distributed along a main level, which can be identified with $Y\times I$ itself, and perhaps several upper levels, which come in two types, depending on whether they correspond to the symplectization of $(Y,\alpha_+)$ or $(Y,\alpha_-)$ (see Figure \ref{coolbuilding}). In further sections, we will refer to these upper levels as \emph{right} or \emph{left}, respectively. 

\begin{figure}[t]\centering
\includegraphics[width=0.45\linewidth]{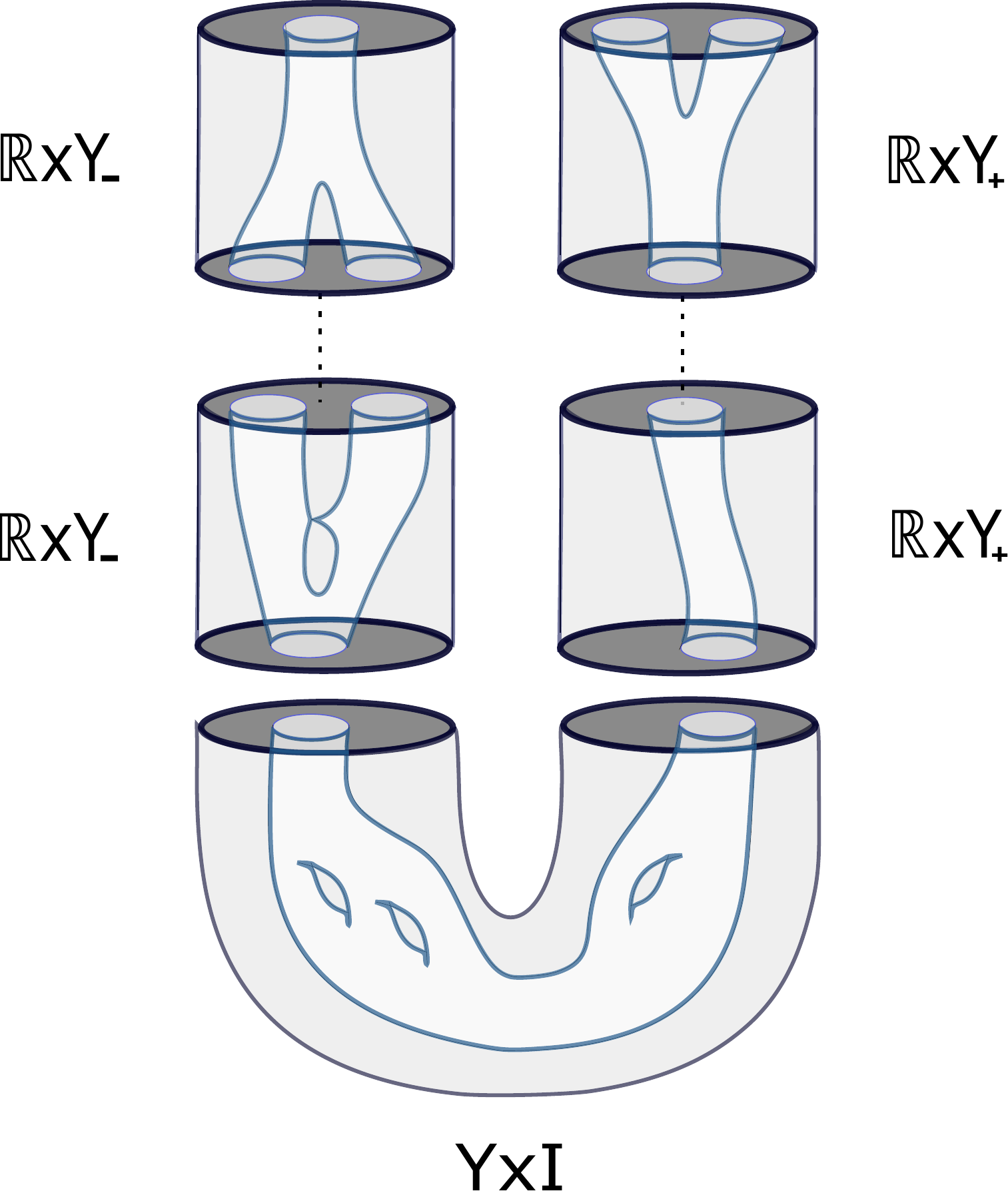}
\caption{\label{coolbuilding} A building in $W_0=\widehat{Y\times I}$, the completion of the cylindrical Liouville semi-filling $Y\times I$.}
\end{figure}

\end{remark}

\subsection{Regularity inside a hypersurface vs. regularity in the symplectization for curves of genus zero}

Consider $u:\dot{S}\rightarrow \mathbb{R}\times M$ a finite energy $J_\epsilon^1$-holomorphic curve, with asymptotics of the form $\gamma_p^l$, with $l\leq N_\epsilon$, and genus $g$. By Proposition \ref{uniqhyp}, $u(\dot{S})\subseteq H$, for some $H \in \mathcal{F}$. If $H$ is positively/negatively asymptotically cylindrical over $H_{p_\pm}$, then the positive/negative asymptotics of $u$ are of the form $\gamma^l_{p_\pm}$, for some simply covered orbit $\gamma\subseteq H_\pm$. Denote by $\Gamma=\Gamma^+\cup \Gamma^-$ the set of punctures of $u$ (where $\Gamma^\pm$ is the set of positive/negative punctures). Define a partition $\Gamma^\pm=\Gamma^\pm_0\cup \Gamma^\pm_1 \cup \Gamma^\pm_2$, where $\Gamma_i^\pm$ denotes the set of positive/negative punctures for which the corresponding asymptotic is of the form $\gamma_p^l$, with $\mbox{ind}_p(H_\epsilon)=i \in \{0,1,2\}$. 

\begin{lemma}\label{regvs}
Assume that $u\subseteq H$ has genus $g=0$, is Fredholm regular inside $H$, and that 
\begin{equation}\label{ATineq2}
\#\Gamma_0^-+\#\Gamma_1^-+\#\Gamma_1^++\#\Gamma_2^-<2
\end{equation}
Then $u$ is Fredholm regular in $\mathbb{R}\times M$.
\end{lemma}

\begin{remark}
Observe that, since all positive/negative asymptotics of $u$ correspond to the same critical point of $H_\epsilon$, one has that $\#\Gamma_i^\pm \neq 0$ implies $\#\Gamma_j^\pm =0$ for $j \neq i$. 
\end{remark}

\begin{proof}
Consider the complex splitting $N_u = N_u^H \oplus u^*N_H$, where $N_u$ is the (generalized) normal bundle to $u$, $N_u^H$ is the normal bundle to $u$ inside $H$, and $N_H$ is the normal bundle to $H$ inside $\mathbb{R}\times M$ (a trivial complex line bundle). Then we get a matrix representation of the normal linearised Cauchy--Riemann operator $\mathbf{D}_u^N: W^{1,2,\delta_0}(N_u)\oplus X_\Gamma \rightarrow L^{2,\delta_0}(\overline{Hom}_\mathbb{C}(\dot{S}, N_u))$, given by 
$$\displaystyle
\mathbf{D}_u^N=\left(\begin{array}{cc} 
\mathbf{D}_{u;T_H}^N & \mathbf{D}_{u;NT}^N  \\
\mathbf{D}_{u;TN}^N& \mathbf{D}_{u;N_H}^N 
\end{array}\right),$$ where 
$$
\mathbf{D}_{u;T_H}^N:W^{1,2,\delta_0}(N_u^H) \oplus X_\Gamma\rightarrow L^{2,\delta_0}(\overline{Hom}_{\mathbb{C}}(\dot{S},N_u^H))
$$
and
$$
\mathbf{D}_{u;N_H}^N:W^{1,2,\delta_0}(u^*N_H) \rightarrow L^{2,\delta_0}(\overline{Hom}_{\mathbb{C}}(\dot{S},u^*N_H))
$$
are Cauchy--Riemann type operators, and the off-diagonal operators are tensorial. Here, $X_\Gamma$ is a finite-dimensional vector space of smooth sections which are supported near infinity and parallel to the Morse--Bott submanifolds containing the asymptotics of $u$. WLOG we can assume that they are tangent to $H$.

The assumption that $u$ is regular inside $H$ is equivalent to the surjectivity of $\mathbf{D}_{u;TH}^N$. Moreover, we claim that $\mathbf{D}_u^N$ is upper-triangular in this matrix representation, i.e.\ that 
$$
\mathbf{D}_{u;TN}^N=0
$$
This, which is basically a generalization of Lemma 3.8 in \cite{Wen1}, can be seen as follows: Take any metric making the splitting $T(\mathbb{R}\times M)\vert_H=T_H \oplus N_H$ orthogonal, and take $\nabla$ to be the associated Levi--Civita connection. This connection is symmetric and preserves this orthogonal splitting. If $\xi \in W^{1,2,\delta_0}(N_u^H)$, take a smooth path $t \rightarrow u_t$ of maps of class $W^{1,p,\delta}(\dot{S},H)$, with image inside $H$, such that $\partial_t\vert_{t=0}u_t=\xi$, $u_0=u$. Since $H$ is holomorphic, $\overline{\partial}u_t \in L^2(\overline{Hom}_{\mathbb{C}}(\dot{S},u_t^*TH))$ takes values in $u_t^*TH$. Then $\mathbf{D}_u^N\xi=\pi_N\nabla_t(\overline{\partial}u_t)\vert_{t=0}$ takes values in $N_u^H \subseteq u^*TH$, since $\nabla$ preserves the splitting, and the claim follows. 

Therefore, to show that $\mathbf{D}_u^N$ is surjective, it suffices to show that $\mathbf{D}_{u;N_H}^N$ is. For this, we can use the automatic transversality criterion in \cite{Wen1}, for the case of a line bundle. We need to check that 
\begin{equation}\label{ATineq}
\mbox{ind}(\mathbf{D}_{u;N_H}^N)>-2+2g+\#\Gamma_{even},
\end{equation}
where we denote by $\Gamma_{even}$ the set of punctures $z$ at which the Conley--Zehnder index of the asymptotic operator of $\mathbf{D}_{u;N_H}^N$ at $z$ is even. Observe that $\Gamma_{even}=\Gamma^+_1 \cup \Gamma^-_1$. Moreover, the operator $\mathbf{D}_{u;N_H}^N$ is asymptotic at each puncture $z$ of $u$ to the normal asymptotic operator associated to the corresponding Reeb orbit $\gamma_{p_z}^{l_z}$ at $z$, whose Conley--Zehnder index is $\mu_N(\gamma_{p_z}^{l_z})=|\mbox{ind}_{p_z}(H_\epsilon)-1|$. By the Riemann--Roch formula, the Fredholm index of $\mathbf{D}_{u;N_H}^N$ is
\begin{equation}\label{indN}
\begin{split}
\mbox{ind}(\mathbf{D}_{u;N_H}^N)=&\chi(\dot{S})+\sum_{z \in \Gamma^+}\mu_N(\gamma_{p_z}^{l_z})-\sum_{z \in \Gamma^-}\mu_N(\gamma_{p_z}^{l_z})\\
=&2-2g-\#\Gamma+\#\Gamma^+_0+\#\Gamma^+_2-\#\Gamma_0^--\#\Gamma_2^-\\
=&2-2g-2\#\Gamma_0^--2\#\Gamma_2^--\#\Gamma^+_1-\#\Gamma_1^-,
\end{split}
\end{equation}
where we used that $\#\Gamma=\sum_{i\in \{0,1,2\}}\#\Gamma^\pm_i$ and $c_1(N_H)=0$ in the natural trivialization. Using that $\#\Gamma_{even}=\#\Gamma_1^++\#\Gamma_1^-$, one sees that the inequality (\ref{ATineq}) can only be satisfied if $g=0$, in which case it is equivalent to (\ref{ATineq2}). 
\end{proof}

\begin{remark} Lemma \ref{regvs} implies that the holomorphic cylinders of Remark \ref{Weinstein} are regular index 1 cylinders in $\mathbb{R}\times M$. This generalizes the situation in \cite{LW}. \end{remark}

\subsection{Obstruction bundles}\label{obbundles}

In this section we will deal with the existence of obstruction bundles over buildings of holomorphic curves in $\mathbb{R}\times M$, to deal with the cases where Fredholm regularity fails. Given such a building $\mathbf{u}$, we want to compute the number of gluings one obtains from $\mathbf{u}$ for a generic perturbation of $J$, in the case where its components fail to be transversely cut-out, but \emph{not too badly}. For this, we require that each component of $\mathbf{u}$ consists of either a regular curve, or a curve for which the dimension of the cokernel of the (normal) linearized Cauchy--Riemann operator is constant as the curve varies in its corresponding moduli space.  We will refer to these curves as \emph{not-too-bad}. If such is the case, one can construct an \emph{obstruction bundle} $\mathcal{O}_{\mathbf{u}}$, and orientable bundle whose fiber is the direct sum of the cokernels of these operators, where the sum varies over the components of $\mathbf{u}$ which are not regular. The base of this bundle is the domain of a \emph{pregluing map}, together with the parameter space keeping track of a fixed perturbation of $J$. Given such a generic (small) perturbation $t\mapsto J_t$ of $J=J_0$ through cyilindrical almost complex structures, the idea is to preglue the components of the $J_0$-holomorphic building $\mathbf{u}$, and impose that the resulting preglued curve be $J_t$-holomorphic. One can view the resulting equation as an obstruction to gluing in the form of a section of the obstruction bundle, whose algebraic count of zeroes is precisely the number of holomorphic gluings of $\mathbf{u}$ one obtains for sufficiently large gluing parameters. This technique, coming from algebraic geometry and Gromov--Witten theory, has been used e.g.\ in \cite{HT1} and \cite{HT2} in the context of ECH (see also \cite{Hut} and \cite{HG} for more gentle introductions).

We sketch the construction in a fairly general case, where curves are immersed and non-nodal. The non-immersed case is slighlty more technical, and for simplicity we will omit it. One needs to replace the normal bundle with the restriction of the tangent bundle of $W$ to the curves in question, including a suitable Teichm\"uller slice parametrizing complex structures on the domain, and divide by the action of the automorphism group. The nodal case just adds and extra $S^1$ as a gluing parameter for each node to the base of the obstruction bundle.

\paragraph*{Setup.} Let $t\mapsto J_t$ be a generic perturbation of $J=J_0$, $t \in [0,1]$, and consider a $J_0$-holomorphic building $\mathbf{u}=(u_1^1,\dots,u_{r_1}^1,\dots, u^N_{1},\dots,u_{r_N}^N)$ with $N\geq 1$ floors, $M \geq N$ components (where $M=\sum_{i=i}^N r_i$) denoted $u_j^i$, consisting of $\mathbb{R}$-translation classes of either regular curves or curves which are not-too-bad. We denote the not-too-bad curves by $u_k$, for $k \in \{1,\dots,r\}$. Assume all of the $u^i_j$ are immersed. We assume that $\mathbf{u}$ contains no trivial cylinder components, since they are regular and glue uniquely, so they do not affect the counts which we will consider. Then there is an orientable obstruction bundle 
$$
\mathcal{O}_{\mathbf{u}}\rightarrow [0,1]\times \prod_{i=1}^{N-1}[R_i,+\infty)\times \prod_{i=1}^N\left(\left(\prod_{j=1}^{r_i} \mathcal{M}^i_j\right)\Big/\mathbb{R}\right):=\mathcal{B}_{\mathbf{u}},
$$
where $R_i\gg 0$ is a gluing parameter for gluing the $i-th$ floor to the $i+1-th$ floor of $\mathbf{u}$, and $\mathcal{M}^i_j$ is the moduli space of $J_0$-holomorphic curves containing $u^i_j$ (without modding out the $\mathbb{R}$-action). The fiber of $\mathcal{O}_\mathbf{u}$ over an element in the base is $\bigoplus_{k=1}^r \mbox{coker}\;\mathbf{D}_k$, where $\mathbf{D}_k$ is the normal linearized Cauchy--Riemann operator of the $k$-th not-too-bad curve $u_k$. 

The base of the obstruction bundle is that of the pregluing map, and the interval $[0,1]$ is the parameter space where $t$ varies. Observe that the case $r=0$ is when every component is regular, and one gets a unique gluing. The case $N=1$ is also considered, in which there is no gluing to do, but the components could be multiply covered, in which case we still wish to count their contributions in the form of a count of zeroes of a section of this bundle.  

\begin{remark}\label{orbifoldcount} \begin{enumerate}[wide, labelwidth=!, labelindent=0pt] $\;$
\item[i)] (Orbibundle) In practice, one should check that the moduli spaces $\mathcal{M}_j^i$ corresponding to not-too-bad curves are orbifolds as worst. The obstruction bundle is then an orientable orbibundle, and the counts of zeroes of a generic section is a weighted, rational count. The way to count is independent of the abstract perturbation scheme for polyfolds (see e.g.\ \cite{HWZ10}, or \cite{Wen3} for a basic exposition on how to count zeroes of sections of orbibundles). In this document, the configurations that we will care about will satisfy this orbifold condition, and existence of obstruction bundles and a way of counting will be enough for our purposes.

\item[ii)] (Tangent space) Under the assumption of i), in practice one should also ask that
\begin{equation}\label{tangentsp}
T_{u_k}\mathcal{M}_k=\ker \mathbf{D}_k,
\end{equation}
where $u_k$ is the $k$-th not-too-bad curve, and $\mathcal{M}_k$ its corresponding moduli space. This condition is automatic for regular curves, but for not-too-bad curves a priori we only get the inclusion $T_{u_k}\mathcal{M}_k\subset\ker \mathbf{D}_k$. In the case of orbifolds, this is to be understood in the orbifold sense, where the tangent space is the quotient of an euclidian space by the linearised action of the isotropy group. The condition is necessary to make sure that the counts of zeroes actually correspond to holomorphic gluings (cf.\ Remark \ref{dimcount}). 

\item[iii)] The fact that $\mathcal{O}_\mathbf{u}$ admits local trivializations can be proven via standard results from non-linear analysis (e.g.\ via the constant rank theorem for Fredholm operators, Corollary 3.1 in \cite{Kai}).

\end{enumerate}
\end{remark}

One has a section $\mathfrak{s}$ of this bundle, the \emph{obstruction section}, defined roughly as follows. Take pregluing data $(t,\mathbf{r}, \mathbf{v}) \in \mathcal{B}_\mathbf{u},$ where
$$\mathbf{r}=(r_1,\dots,r_{N-1}),$$
$$\mathbf{v}=(v^1_1,\dots,v^1_{r_1},\dots,v^N_1,\dots,v^N_{r_N})$$
We identify the $\mathbb{R}$-translation classes $v^i_j$ with representatives, and we consider the components $v^i_{r_i}$ as fixed $\mathbb{R}$-level. Consider a tuple $$\Psi=(\psi^1_1,\dots,\psi_{r_1}^1,\dots, \psi^N_1,\dots,\psi_{r_N}^N),$$ where $\psi^i_j$ is a section of the normal bundle of $v^i_j$ belonging to a suitable Hilbert completion of the space of smooth sections. 

One constructs a (only \emph{approximately} $J_t$-holomorphic) preglued curve $\oplus_{(t,\mathbf{r},\Psi)}\mathbf{v}$ out of the $J_0$-building $\mathbf{v}$, with the property that it converges to $\mathbf{v}$ as every component of $\mathbf{r}$ approaches $+\infty$. This is a standard construction, and is done by exponentiating $\psi^i_j$ along $v^i_j$, and by using suitable smooth cut-off functions $\beta^i_j$, which equal $1$ in the interior of $v^i_j$ and decay towards its cylindrical ends, translated in a way determined by the gluing parameters $\mathbf{r}$. 

If we impose that $\oplus_{(t,\mathbf{r},\Psi)}\mathbf{v}$ be $J_t$-holomorphic, we get an equation of the form 
$$
\bigoplus_{i=1}^N\bigoplus_{j=1}^{r_i} \beta^i_j \Theta^i_j(t,\mathbf{r},\Psi)=0
,$$ 
where
$$
\Theta^i_j(t,\mathbf{r},\Psi)= \mathbf{D}^i_j\psi^i_j+\mathcal{F}^i_j(t,\mathbf{r},\Psi)
$$
Here, $\mathbf{D}^i_j$ is the linearized Cauchy--Riemann operator of $v^i_j$ for $t=0$, and the second summand involves extra terms arising from the patching construction, mostly non-linear and depending only on the $\psi^m_n$'s for which $v^i_j$ is adjacent to $v^m_n$. For sufficiently small $t$, and fixed sufficiently large $R_i$, there exists a unique $\Psi=\Psi_{t,\mathbf{v}}$ such that 
$$
\psi^i_j \perp \ker\mathbf{D}^i_j, \mbox{ for every } i,j. 
$$
$$
\Theta^i_j(t,\mathbf{r},\Psi)=0, \mbox{ if } u^i_j \mbox{ (or equivalently } v^i_j \mbox{) is regular\ }.
$$
$$
\Theta^i_j(t,\mathbf{r},\Psi) \perp \mbox{im}\;\mathbf{D}^i_j, \mbox{ if } u^i_j \mbox{ (or equivalently } v^i_j \mbox{) is not-too-bad\ }.
$$
Then one can define 
$$
\mathfrak{s}(t,\mathbf{r},\mathbf{v})=\bigoplus_{k=1}^r \pi_k  \Theta_k(t,\mathbf{r},\Psi_{t,\mathbf{v}})=\bigoplus_{k=1}^r \pi_k  \mathcal{F}_k(t,\mathbf{r},\Psi_\mathbf{t,v}),
$$
where $\pi_k$ is the orthogonal projection to $\mbox{coker}\; \mathbf{D}_k$. Therefore, there is a 1-1 correspondence between the zeroes of $\mathfrak{s}(t,\cdot,\mathbf{v})$ and the $J_t$-holomorphic gluings of $\mathbf{v}$. Moreover, $\mathfrak{s}(t,\cdot,\mathbf{v})$ is transverse to the zero section for generic $J_t$. 

\begin{remark}\label{dimcount}

Under the assumptions of Remark \ref{orbifoldcount}, one easily computes that
\begin{equation}
\begin{split}
\dim \mathcal{B}_u&=\mbox{virt-dim}(\mathbf{u}) + \mbox{rank} \;\mathcal{O}_\mathbf{u}
\end{split}
\end{equation}

We conclude that if $\mbox{virt-dim}(\mathbf{u})=1$, then generically, for \emph{fixed} small $t$ and fixed large gluing parameters, there will finitely many $J_t$-gluings.  

\end{remark}

\subparagraph*{Existence of obstruction bundles.}\label{OBB} In the case where the hypothesis of the Lemma \ref{regvs} are not satisfied for a possibly multiply-covered curve $u:(\dot{S},j)\rightarrow H\subseteq W=\mathbb{R}\times M$, we wish to have a criterion for the existence of an obstruction bundle for building configurations containing $u$. We will assume that $u$ is is not-too-bad in $H$, which includes the case where $u$ is regular. We prove that $u$ is not-too-bad in $\mathbb{R}\times M$, which provides obstruction bundles in $\mathbb{R}\times M$.

As in the proof of Lemma \ref{regvs}, consider the splitting 
\begin{equation}
\label{splitop}
\displaystyle
\mathbf{D}_u^N=\left(\begin{array}{cc} 
\mathbf{D}_{u;T_H}^N & \mathbf{D}_{u;NT}^N  \\
0 & \mathbf{D}_{u;N_H}^N 
\end{array}\right)
\end{equation}

If the hypothesis of Lemma \ref{regvs} fail, then we have 
$$
\mbox{ind}\;\mathbf{D}_{u;N_H}^N\leq 0
$$
Indeed, this follows from equation (\ref{indN}), if we assume that $g\geq 1$ or $\#\Gamma_0^-+\#\Gamma_1^-+\#\Gamma_1^++\#\Gamma_2^-\geq 2$. By \cite{Wen1}, we get a bound
\begin{equation}\label{boundker}
\dim \ker \mathbf{D}_{u;N_H}^N \leq K(c_1(N_H,u),\#\Gamma_{even}),
\end{equation}
where 
$$
2c_1(N_H,u)=\mbox{ind}\;\mathbf{D}_{u;N_H}^N-2+2g+\#\Gamma_{even}=-2(\#\Gamma_0^-+\#\Gamma_2^-)\leq 0
$$
and
$$
K(c,G)=\min\{k+l|k,l \in \mathbb{Z}^{\geq 0}, k\leq G, 2k+l>2c, \; l \mbox{ even}\}
$$

In particular, if either $\#\Gamma_0^-$ or $\#\Gamma_2^-$ are non-zero (which is the case only when $u$ lies in the cylindrical hypersurface over $min$ or $max$), then $\mathbf{D}_{u;N_H}^N$ is injective.  In the case where $c_1(N_H,u)=0$, we obtain 
$$
\dim \ker \mathbf{D}_{u;N_H}^N\leq 2
$$
If $H$ is non-cylindrical, then $\dim \ker \mathbf{D}_{u;N_H}^N \geq 1$, since $\partial_a$ is then normal to $H$, and the almost complex structure is $\mathbb{R}$-invariant. If $H$ projects to a flow-line joining the maximum to the minimum, all nearby index $2$ flow-lines are obtained by a push-off in the $\Sigma$-direction of a normal section $\nu$ which decays asymptotically, and this corresponds to holomorphic push-offs of $u$ in nearby hypersurfaces. We conclude that 
$$
\dim \ker \mathbf{D}_{u;N_H}^N=2
$$
In the non-generic case of an index $1$ flow-line, the latter section is not there, and one has $\dim \ker \mathbf{D}_{u;N_H}^N=1,$ spanned by the $\mathbb{R}$-direction. Using that
$$
\ker \mathbf{D}_u^N = \ker \mathbf{D}_{u;T_H}^N \oplus \langle \partial_a,\eta \rangle,
$$
in the generic case of an index $2$ flow-line, and 
$$
\ker \mathbf{D}_u^N = \ker \mathbf{D}_{u;T_H}^N \oplus \langle \partial_a \rangle,
$$
in the non-generic case of an index $1$ flow-line, and that the index is only dependent on the moduli space, we conclude the same for the cokernel. This finishes the non-cylindrical case.

The last case is when $H$ is cylindrical over a hyperbolic critical point, which doesn't follow from automatic transversality. We show that the operator $\mathbf{D}_{u;N_H}^N$ is injective by a perturbation argument, when $H$ is cylindrical (not necessarily over a hyperbolic point), as follows. For $\epsilon=0$, $\mathbb{R}\times M_P^\pm$ is foliated by cylindrical holomorphic hypersurfaces $\mathbb{R}\times Y \times \{p\}$ for \emph{any} $p\in M_P^\pm$. This implies that

\[
\displaystyle
\mathbf{D}_u^N=\mathbf{D}_{u,\epsilon}^N=\left(\begin{array}{cc} 
\mathbf{D}_{u;T_H}^N & \mathbf{D}_{u;NT}^N  \\
0 & \mathbf{D}_{u;N_H}^N 
\end{array}\right)\stackrel{\epsilon\rightarrow 0}{\longrightarrow}
\left(\begin{array}{cc} 
\mathbf{D}_{u;T_H}^N & 0  \\
0 & \overline{\partial} 
\end{array}\right)
\]

The operator $\overline{\partial}$ is injective, since the elements in its kernel are holomorphic sections which decay at the punctures. Since injectivity is an open condition in the usual Fredholm operator topology, it follows that $\mathbf{D}_{u;N_H}^N$ is injective for sufficiently small $\epsilon>0$.
Therefore 
$$ 
\ker \mathbf{D}_u^N=\ker \mathbf{D}_{u;T_H}^N,
$$
which depends only on the moduli containing $u$, by assumption. This finishes all cases. 

However, the small $\epsilon$ needs to be depends a priori on the curve $u$. While this is perhaps just a technicality, as long as we consider curves $u$ with Morse--Bott asymptotics, of which the positive have bounded total action (for a fixed bound), we can assume that this operator is injective for such family of curves, since there will be finite such moduli spaces. 

Observe that in the case of a regular cylinder inside a cylindrical hypersurface, this implies that it is regular in $\mathbb{R}\times M$. This includes the case of a cylinder lying in a cylindrical hypersurface corresponding to a hyperbolic critical point, which is not covered by Lemma \ref{regvs}.   

We have proved the following:

\begin{prop}[cylinders]\label{regcylcyls}
Every cylinder with two positive ends over a hyperbolic critical point, and regular in a non-cylindrical hypersurface, has an obstruction bundle of rank 1. We can take $\epsilon>0$ sufficiently small so that cylinders which are regular inside their corresponding hypersurface are also regular in $\mathbb{R}\times M$, for every other case. For cylinders lying in cylindrical hypersurfaces corresponding to a hyperbolic critical point, we can ensure their regularity as long as we consider fixed action bounds on the positive asymptotics, or finite families of moduli spaces. 
\end{prop} 

More generally,

\begin{prop}\label{exisob} Assume that $u\subseteq H$ is not-too-bad inside $H$, and the rest of the hypothesis of Lemma \ref{regvs} fail (i.e.\ either $g>0$ or inequality (\ref{ATineq2}) fails). If $u$ does not lie in a cylindrical hypersurface corresponding to a hyperbolic critical point, then $u$ is not-too-bad in $\mathbb{R}\times M$. Moreover, for every $T>0$ action threshold, we can choose an $\epsilon>0$ sufficiently small, such that every curve $u$ which lies in such a hypersurface and is not-too-bad inside of it, and such the total action of its positive asymptotics is bounded by $T$, is not-too-bad in $M$.

In all above cases, this implies that there exists an obstruction bundle $\mathcal{O}_u$ for gluing $u$ to any building configuration which contains it.

\end{prop}

\begin{remark} It is not hard to show that, for small $\epsilon>0$, the rank of $\mathcal{O}_u$ is
\begin{equation}\label{rankO}
\mbox{rank}\;\mathcal{O}_u=\mbox{rank}\;\mathcal{O}^H_u -\mbox{ind}\;\mathbf{D}_{u;N_H}^N+\dim \ker \mathbf{D}_{u;N_H}^N, 
\end{equation}
where $\mathcal{O}^H_u$ is the obstruction bundle of $u$ inside $H$, and the second term is given by formula (\ref{indN}). Recall that in the case for non-cylindrical hypersurfaces $H$, we have $\dim \ker \mathbf{D}_{u;N_H}^N\in \{1,2\}$ depending on whether the corresponding flow line is index $1$ or $2$; and for cylindrical hypersurfaces, $\dim \ker \mathbf{D}_{u;N_H}^N=0$. One can show also that
\begin{equation}\label{cokeq}
\mbox{coker}\;\mathbf{D}_u^N=\mbox{coker}\; \mathbf{D}_{u;T_H}^N \oplus \mbox{coker}\; \mathbf{D}_{u;N_H}^N
\end{equation}
\end{remark}

\begin{remark}
Given a leafwise not-too-bad curve $u$, if we assume that condition (\ref{tangentsp}) holds for the leafwise moduli space and the operator $\mathbf{D}^N_{u;T_H}$, then our computations of the kernel of the total operator $\mathbf{D}_u^N$ implies that it also holds in the ambient manifold $\mathbb{R}\times M$. 
\end{remark}

This finishes the general construction. From now on, we discuss a particular subclass of examples, for which we obtain the results from the introduction. 

\section{Non-fillable 5-dimensional model with no Giroux torsion}\label{5dmodel}

For this chapter, we fix $Y=ST^*X$ to be the unit cotangent bundle of a hyperbolic surface $X$ with respect to a choice of hyperbolic metric, and $\pi: Y \rightarrow X$ the natural projection. The 1-form $\alpha_-$ is the standard Liouville form, and $\alpha_+$ is a prequantization contact form. This means that $\alpha_+$ is a connection form with curvature $d\alpha_+=\pi^* \omega$, where $\omega$ is a symplectic form on $X$ representing $c_1(Y)=-2+2g(X)>0$ in $H^2(X;\mathbb{R})\simeq\mathbb{R}$ when $Y$ is viewed as a complex line bundle, i.e.\ area$(\omega)=\int_X \omega = -\chi(X)$. This example of $Y \times I$ was originally constructed in \cite{McD}. We consider the family of $5$-dimensional contact manifolds $M=Y \times \Sigma$ constructed in previous sections.

We will dig into the SFT of this class of examples, and derive our results. We shall investigate whether these examples have 1-torsion (for any $k$, not just $k=1,2$ in which case we already know they do). In Section \ref{proooof}, we will classify all possible building configurations that can contribute to 1-torsion in the whole SFT algebra, and prove Theorem \ref{nocob1} in Section \ref{proofofnocob}. 
  
\subsection{Curves on symplectization of prequantization spaces: Existence and uniqueness}
\label{prequantcurves}

Consider $(Y,\xi=\ker e^H \lambda)$ a prequantization space over an integral symplectic base $(X,\omega)$, where $H$ is a Hamiltonian on $X$, and $d\lambda=\pi^*\omega$. It is a reasonably standard construction that any choice of $\omega$-compatible complex structure on $X$ naturally lifts to a compatible $J$ on the symplectization of $Y$, for which there exists a finite energy foliation $\mathcal{F}_H$ by $J$-holomorphic cylinders. These project to $X$ as flow lines of $H$, and their asymptotics correspond to circle fibers over the critical points of $H$ (see e.g.\ \cite{Sie, Mo2}). 

\begin{lemma}[Uniqueness]
\label{uniqpre}
Let $X$ be a closed surface of genus $g\geq 1$, and let $\omega$ be an integral area form on $X$. Let $\pi:Y\rightarrow X$ denote the prequantization space over $X$, with a connection (contact) form $\lambda$ whose curvature form satisfies $d\lambda=\pi^*\omega$. Then, for any choice of $\omega$-compatible complex structure on $X$ lifting to a compatible almost complex structure $J$ on $\mathbb{R}\times Y$,  any action threshold $T>0$, and any Morse Hamiltonian $H:X\rightarrow \mathbb{R}$, we can find sufficiently small $\epsilon>0$, such that any genus  $g^\prime<g$ holomorphic curve on $(\mathbb{R}\times Y,e^{\epsilon H}\lambda)$ whose positive asymptotics all correspond to critical points of $H$, and have total action bounded by $T$, is a multiple cover of a flow line cylinder in the finite energy foliation $\mathcal{F}_H$ induced by $H$.
\end{lemma}

\begin{proof} 
In the degenerate case $H=0$, the projection $\mathbb{R} \times Y \rightarrow X$ is holomorphic, and every Reeb orbit in $Y$ is a multiple of the $S^1$-fiber. So, a curve as in the statement induces a holomorphic map into $X$, defined on a closed curve of genus $g^\prime$. By holomorphicity, it has non-negative degree, and has zero degree if and only if it is constant. If it has positive degree, by Poincar\'e duality it induces an injection in cohomology. But this cannot happen if $g^\prime<g$. We conclude that it is a cover of a the trivial cylinder. 

In the non-degenerate case where $H$ is small but non-identically zero Morse perturbation, one can use holomorphic cascades to reduce to the degenerate case, or alternatively the results from \cite{MS} (see \cite{Mo2}). 
\end{proof}

\subsection{Some remarks and useful facts}

Recall that $(Y=ST^*X, \alpha_-)$ is the unit cotangent bundle of a hyperbolic surface, and $\alpha_-$ is the standard Liouville form. We shall need the following facts about closed $\alpha_-$-orbits (e.g.\ \cite{CL09,Luo,Sch}):
\begin{itemize}
\item They project to $X$ as closed hyperbolic geodesics.
\item They are non-degenerate.
\item Their Conley--Zehnder index vanish in a natural trivialization.
\end{itemize}

Take $J_0$ a cylindrical $d\alpha$-compatible almost complex structure in $W_0=\widehat{Y\times I}$. Recall that in Section \ref{compJmodelB} we have defined an almost complex structure $J_\epsilon^1$ in $\mathbb{R}\times M$, which is compatible with the SHS $\mathcal{H}_\epsilon^1$, and for which we have a foliation by holomorphic hypersurfaces $\mathcal{F}$. By a generic perturbation of $J_0$ (and hence of $J_\epsilon^1$ along the hypersurfaces in $\mathcal{F}$), we can assume that:

\begin{itemize}
\item[\textbf{(A)}] Every somewhere injective holomorphic curve $u$ in the completion $W_0=\widehat{Y\times I}$ which intersects the main level is regular, and therefore satisfies ind$_{W_0}(u)\geq 0$ (where this denotes the index computed in $W_0$). 
\end{itemize}

Using that $c_1(\xi_\pm)=0,$ from Remark \ref{c1rk}, we obtain
\begin{itemize} 
\item[\textbf{(B)}] $c_1(\xi_\epsilon^1)=0$.
\end{itemize} 

Observe that, since closed Reeb orbits for $\alpha_-$ project down to closed hyperbolic geodesics, they are non-contractible, and Reeb orbits for $\alpha_+$ project to points. This implies the following fact:

\begin{itemize}
\item[\textbf{(C)}] There are no holomorphic cylinders in our model crossing sides of the dividing set.
\end{itemize}

This does not happen for the $3$-dimensional cases in \cite{LW}, where the hypersurfaces are cylinders.

Observe that we are allowed to take a \emph{generic} almost complex structure $J_-$ on $(\mathbb{R}\times Y_-,d(e^a\alpha_-))$ and use it in the construction of our model. Since the index of every cylinder is necessarily zero, every multiply-covered cylinder is unbranched, and multple covers of trivial cylinders are again trivial, we obtain:

\begin{itemize}
\item[\textbf{(D)}] Every holomorphic cylinder in $(\mathbb{R}\times Y, d(e^a\alpha_-))$ is necessarily trivial. And every holomorphic curve of index zero is necessarily a trivial cylinder.
\end{itemize}

Recall that for a given free homotopy class of simple closed curves in a hyperbolic surface $X$, there is a unique geodesic representative. This implies that if two Reeb orbits in $(Y,\alpha_-)$ are joined by a holomorphic cylinder in $W_0$ with two positive ends (which by (\textbf{D}) is necessarily contained in the main level $Y\times I$), they correspond to the same geodesic, but with different orientation. In particular, they are \emph{disjoint} Reeb orbits. We conclude:

\begin{itemize}
\item[\textbf{(E)}] In $W_0$, every holomorphic cylinder with two positive ends on a left upper level joins disjoint Reeb orbits.
\end{itemize}

\subsection{Investigating 1-torsion}\label{proooof}

In this section, we will study the portion of the SFT differential of the contact manifold $(M=ST^*X\times \Sigma,\xi)$ which has the potential to yield $1$-torsion, where we denote by $\xi$ the isotopy class of contact structures on $M$ that we have defined in previous sections. Recall from the introduction that we have a power series expansion $\mathbf{D}_{SFT}=\sum_{k\geq 1}D_k \hbar^{k-1}$ of said differential, where $D_k$ is a differential operator of order $\leq k$. This operator is defined by counting holomorphic curves in $\mathbb{R}\times M$ with $|\Gamma^+|+g=k$, where $\Gamma^+$ is the set of positive punctures, and $g$ is the genus. Here, we consider the untwisted version of the SFT, where we do not keep track of homology classes. We will compute the projection of the operators $D_1$ and $D_2$ to the base field $\mathbb{R}$, but considering their actions on Reeb orbits only up to a large action threshold. In other words, we forget the $q$-variables from these operators, so that we only consider curves with no negative ends. While the projection of $D_1$ to $\mathbb{R}$ vanishes identically (there are no holomorphic disks), the computation for $D_2$ is rather involved. The result is rather unexpected: among precisely $35$ possibilities, there is basically only one way of obtaining $1$-torsion. While we cannot prove rigorously that $1$-torsion indeed arises, in the next section we provide a heuristic argument as to why we expect this to be true. While we originally expected that the classification of contributions would show that $1$-torsion does not arise, we will use our knowledge of cylinder configurations to prove Theorem \ref{nocob1}.

\paragraph*{Setup.} Denote by $H_X$ a choice of a Morse function satisfying the Morse--Smale condition on $X$. Choose a (non-generic) cylindrical almost complex structure $J_+$ on the symplectization of $(Y,\alpha_+=e^{\epsilon^\prime H_X}\lambda)$, coming from a lift of a complex structure on $X$, inducing a foliation $\mathcal{F}_{X}$ of $\mathbb{R}\times Y$ by $H_X$-flow-line cylinders (as described in Section \ref{prequantcurves}). Here, choose $\epsilon^\prime>0$ small enough so that Lemma \ref{uniqpre} holds for $T=T_{\epsilon^\prime}$. Denote $\gamma_{p;q}$ the simply covered Reeb orbit in $M$ corresponding to $(p,q)\in \mbox{crit}(h_+)\times \mbox{crit}(H_X)$. Given a holomorphic curve (or building) $u$ with asymptotics corresponding to critical points of $h_\pm$, and then lying in a hypersurface of the foliation $\mathcal{F}$, we will view it as a punctured curve in $W_0$, when convenient. We will denote by $\mbox{ind}_M(u)$ (and $\mbox{ind}_{W_0}(u)$) its Fredholm index when viewed as a curve in $\mathbb{R}\times M$ (resp. in $W_0$). By \textbf{(B)}, we have $$\mbox{ind}_M(u)=(\mu^\tau_{CZ})_M(u)$$$$ \mbox{ind}_{W_0}(u)=-\chi(u)+(\mu^\tau_{CZ})_{W_0}(u)$$ By Proposition \ref{indexthm}, given $l\leq N_\epsilon$, then for orbits $\gamma^l_{p;q}$ for hyperbolic $p$, we have $(\mu^\tau_{CZ})_M(\gamma^l_{p;q})= (\mu^\tau_{CZ})_{W_0}(\gamma^l_{p;q})=\mbox{ind}(q)-1$; and for orbits $\gamma^l_{p;q}$, for $p$ the maximum or minimum, we have $(\mu^\tau_{CZ})_M(\gamma^l_{p;q})= (\mu^\tau_{CZ})_{W_0}(\gamma^l_{p;q})+1=\mbox{ind}(q)$. 

For $k\leq 2$, we already know that our model has $1$-torsion, so we set $k\geq 3$. Denote by $\overline{\mathcal{M}}^1_{g,r}(J_\epsilon^1;T_\epsilon)$ the moduli space of all translation classes of index $1$ connected $J_\epsilon^1$-holomorphic buildings in $\mathbb{R}\times M$ with arithmetic genus $g$, no negative ends, and $r$ positive ends approaching orbits whose periods add up to less than $T_\epsilon$. Elements in this moduli are the buildings that potentially glue (after perturbation) to curves which might contribute to the SFT differential. We will prove the existence of a choice of coherent orientations such that elements in $\overline{\mathcal{M}}^1_{g,r}(J_\epsilon^1;T_\epsilon)$, whenever $g+r\leq 2$, cancel in pairs (after introducing all SHS-to-contact, non-degenerate, and a generic perturbation), \emph{except for a single sporadic building configuration}. This unusual configuration consists of a single-level punctured torus, and we will discuss it in the next section. This means that the only way of obtaining 1-torsion is to differentiate the positive asymptotic of one of these configurations. In particular, configurations corresponding to cylinders come in cancelling pairs, and this is what we will use to show that there is no Giroux torsion in our model. 

\paragraph*{Classification.} If $\mathbf{u}$ is an element of $\overline{\mathcal{M}}^1_{g,r}(J_\epsilon^1;T_\epsilon)$ with $g+r\leq 2$, since $r\geq 1$ by exactness, then it can only glue to one of the following: 

\begin{itemize}
\item \textbf{Case 1.} A plane with one positive end.

\item \textbf{Case 2.} A cylinder with two positive ends.

\item \textbf{Case 3.} A 1-punctured torus with one positive end. 

\end{itemize}

Case 1 is ruled out, since Reeb orbits in $Y$ are non-contractible (for both the prequantization form and the standard Liouville form), as are the Reeb orbits in $Y \times \mathcal{U}$. The other orbits are ruled out by the definition of $T_\epsilon$ and $\overline{\mathcal{M}}^1_{g,r}(J_\epsilon^1;T_\epsilon)$. This means that there is no $0$-torsion, so this already shows our model is tight (we already knew this from \cite{Mo2}, since both $(Y,\alpha_\pm)$ are hypertight).

\vspace{0.5cm}

So we deal with the other two cases. Case 2 is the only one which is relevant for Theorem \ref{nocob1} and Corollary \ref{corGT}, so we will not include details for Case 3, which are rather involved. Using that $k\geq 3$ and the action restriction, as in \cite{LW, Mo2}, one easily checks that the Reeb orbits in the asymptotics of every floor of $\mathbf{u}$ are of the form $\gamma_{p;q}^l$, for some $l\leq N_\epsilon$. We can then appeal to Proposition \ref{uniqhyp}, which yields the fact that each component of $\mathbf{u}$ lies in one of the holomorphic hypersurfaces of the foliation, so that we may view them as a building inside $W_0=\widehat{Y\times I}$. 

\vspace{0.5cm} 

\textbf{Case 2. A cylinder with two positive ends.} 

\begin{figure}[t]\centering
\includegraphics[width=0.15\linewidth]{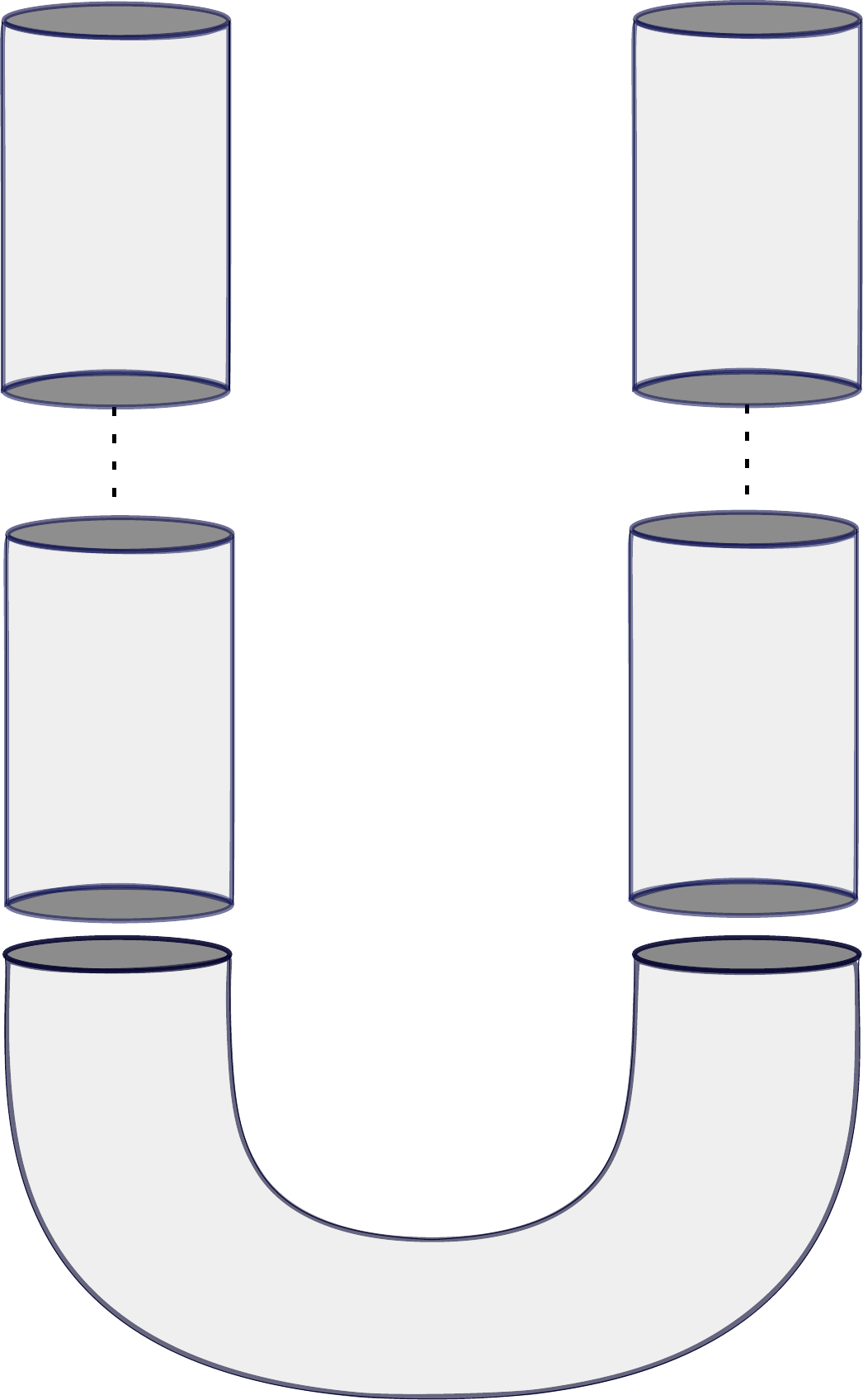}
\caption{\label{poss} The only possible a priori combinatorics for $u$  in the case $g=0$, $r=2$.}
\end{figure}

Since $\mathbf{u}$ glues to a cylinder with no negative ends, and since the Reeb orbits $\gamma_{p;q}$ are non-contractible then one can show that $\mathbf{u}$ can only consist of a cylinder with two positive ends in the bottom floor, together with two chains of cylinders on top of its ends (see Figure \ref{poss}). 

By \textbf{(C)}, no component of $\mathbf{u}$, all cylinders, can have asymptotics in different sides. Also, recall that we may assume that $\mathbf{u}$ is a stable building, so that it does not have levels consisting only of trivial cylinders. 
   
We see that $\mathbf{u}$ cannot consist solely of right upper level components (corresponding to the prequantization space), since its bottom level would be a cylinder with two positive ends and then cannot be a cover of a flow line cylinder, which contradicts Lemma \ref{uniqpre}. We also see that it cannot consist solely of left upper level components, since its bottom level cannot be a trivial cylinder (it has no negative ends), which contradicts observation \textbf{(D)}.

Then $\mathbf{u}$ has a nontrivial component in the main level, which we call $u_0$.  

\vspace{0.5cm}

\begin{figure}[t]\centering
\includegraphics[width=0.50\linewidth]{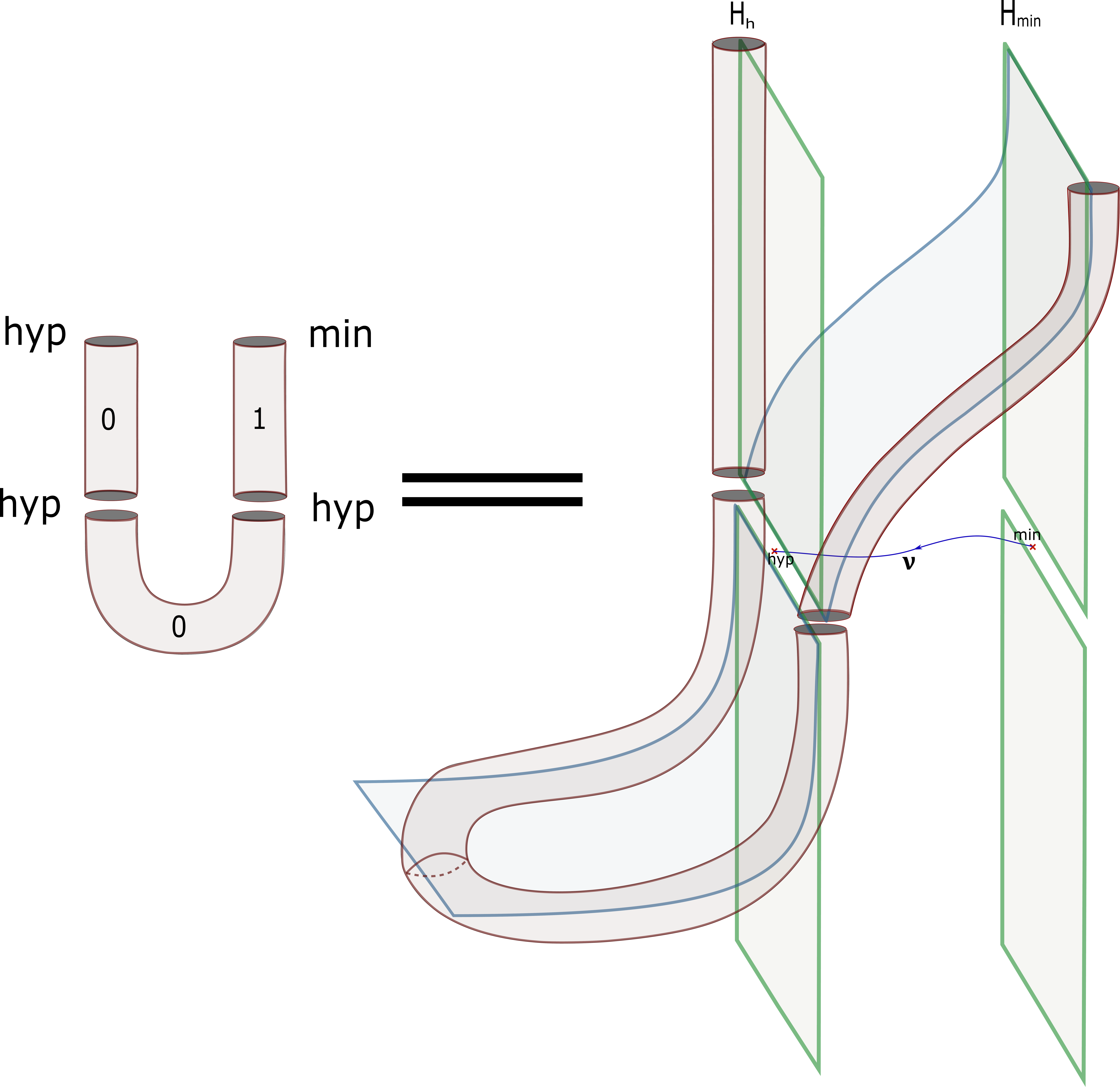}
\caption{\label{excurve}A possible configuration for the index $1$ building $u$ for the case $g=0$, $r=2$. }
\end{figure}

\textbf{Case 2.A. Both asymptotics of $u_0$ lie on a left level.} Label by $\bold{hyp}$ or by $\bold{min}$ the Reeb orbits appearing in $\mathbf{u}$, according to whether they lie over a hyperbolic point, or the minimum. Since $u_0$ lies in a hypersurface, its two positive asymptotics necessarily have the same label. Thus, for each string, the associated ordered sequence of labels can only look like $(\bold{hyp},\dots,\bold{hyp},\bold{min},\dots,\bold{min})$. The number of $\bold{hyp}$'s or $\bold{min}$'s may be zero, but not both (see Figure \ref{excurve}). Observation \textbf{(D)} implies that all the upper components correspond to trivial cylinders (under Remark \ref{Weinstein}). We then see that the only possibility for $\mathbf{u}$ is the one depicted in Figure \ref{excurve}, since all others will have index different from $1$.

Observe that $u_0$ has only hyperbolic asymptotics, and, if it is multiply covered, it is necessarily unbranched. Since $J_0$ is generic, Theorem \ref{PM} implies that $u_0$ is regular in $W_0$, and the same is true for its upper components, which are trivial. By Proposition \ref{regcylcyls}, we obtain an obstruction bundle of rank $1$ for this configuration. 

This configuration has an ``evil twin'', obtained by replacing the index $1$ upper component, which lies over an index $1$ Morse flow-line $\gamma$ connecting $hyp$ to $min$, with the index $1$ cylinder lying over the flow-line $\overline{\gamma}$ (the unique other flow-line connecting $hyp$ to $min$). If we take our coherent orientation to be compatible with the Morse orientation, the algebraic counts of zeroes of the obstruction sections associated to the twin configurations will cancel out after introducing a generic perturbation. One may check that the asssociated obstruction bundles satisfy the properties of Remark \ref{orbifoldcount}.

\vspace{0.5cm}

\textbf{Case 2.B. Both asymptotics of $u_0$ lie on a right level.}

Every component in the upper levels of $\mathbf{u}$ is either contained in a holomorphic hypersurface $H_a^\nu$ corresponding to an index 1 Morse flow line $\nu$ connecting a hyperbolic point $hyp \in \Sigma_+$ to the maximum $max \in \Sigma_+$, or is contained in either a hypersurface $\mathbb{R}\times H_{hyp}$ or $\mathbb{R}\times H_{max}$ lying over $hyp$, or $max$. Moreover, Lemma \ref{uniqpre} implies that, when $\mathbf{u}$ is viewed as lying inside $W_0$, they are (necessarily unbranched) covers of flow-line cylinders. Denote by $q_1,q_2 \in \mbox{crit}(H_X)\subseteq X$ the critical points corresponding to the two positive asymptotics of $u_0$. Label by $(\bold{hyp};\mbox{ind}_q(H_X))$ or by $(\bold{max};\mbox{ind}_q(H_X))$ the Reeb orbits appearing in $u$, according to whether they lie over a hyperbolic point, or the maximum, and where $q \in \crit(H_X)$ is the corresponding critical point. Again, since $u_0$ lies in a hypersurface, the first component of the labels of its two positive asymptotics necessarily agree.

Observe that, if $v$ denotes a right upper level component of $u$, one has 
$$1=\mbox{ind}_M(u)=\mbox{ind}_M(u_0)+\sum_{v}\mbox{ind}_M(v)\geq \mbox{ind}_M(v)\geq \mbox{ind}_{W_0}(v)\geq 0$$
Denote by $u_0^\prime$ the somewhere injective curve underlying $u_0$, which, since Reeb orbits are non-contractible, is a cylinder over which $u_0$ is unbranched, and satisfies $\mbox{ind}_{W_0}(u_0^\prime)=\mbox{ind}_{W_0}(u_0)$. By \textbf{(A)}, we have that 
$$
0\leq \mbox{ind}_{W_0}(u_0^\prime)=\mbox{ind}_{W_0}(u_0)=\mbox{ind}_p(H_X)-1+\mbox{ind}_q(H_X)-1\leq\mbox{ind}_M(u)=1,
$$
so that 
$$3\geq \mbox{ind}_p(H_X)+\mbox{ind}_q(H_X)\geq 2$$

If the two labels of the positive ends of $u_0$ are $\bold{max}$, then its index in $M$ is 
$$
\mbox{ind}_M(u_0)= \mbox{ind}_p(H_X)+\mbox{ind}_q(H_X)\geq 2
$$
Since $1=\mbox{ind}_M(u) \geq \mbox{ind}_M(u_0)\geq 2$, we get a contradiction. Then both asymptotics of $u_0$ have $\bold{hyp}$ as the first component of their label, which in particular implies $\mbox{ind}_{W_0}(u_0)=\mbox{ind}_M(u_0)$. 

\begin{figure}[t]\centering
\includegraphics[width=0.99\linewidth]{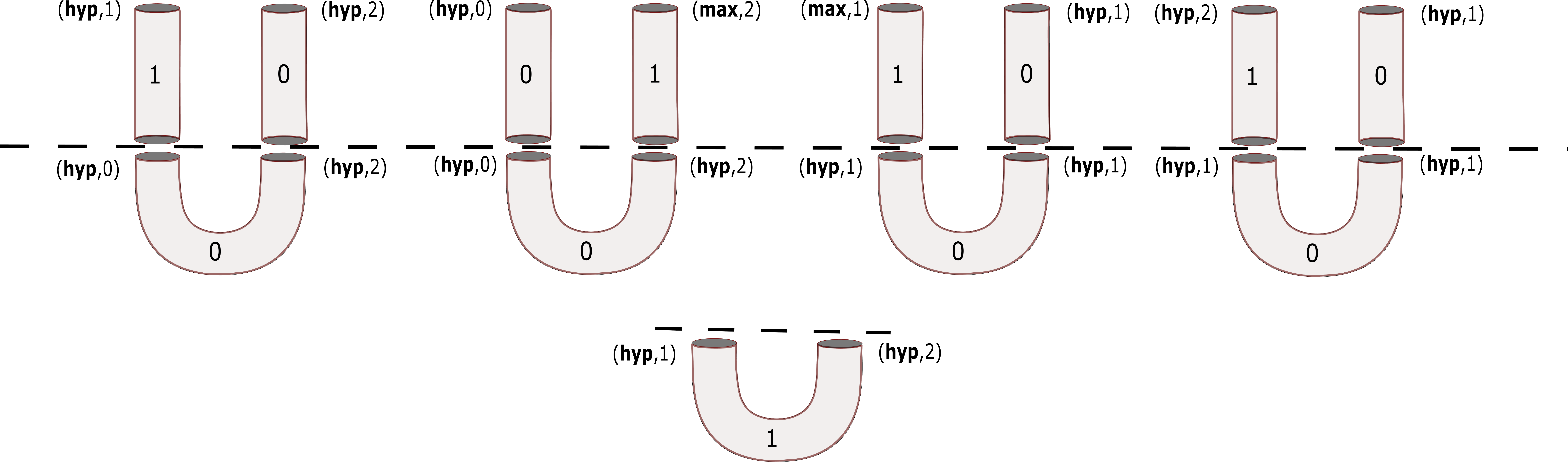}
\caption{\label{cylconf} All 5 possible configurations with positive ends on a right upper level (up to obvious symmetries). The dotted lines separates the main level from the upper levels. The upper level components correspond to unbranched covers of flow-line cylinders over $X$, which are regular in $W_0$.}
\end{figure}

Assume $\mbox{ind}_p(H_X)+\mbox{ind}_q(H_X)=2$, so that $u_0$ has index zero in $W_0$. In the case that $(\mbox{ind}_p(H_X),\mbox{ind}_q(H_X))=(0,2)$, the automatic transversality criterion in \cite{Wen1} implies that it is regular in $W_0$. If $(\mbox{ind}_p(H_X),\mbox{ind}_q(H_X))=(1,1)$, then we use Theorem \ref{PM} to conclude again that $u_0$ is regular in $W_0$. The only possibilities for $\mathbf{u}$ are shown in Figure \ref{cylconf}, all of them have an ``evil twin'' in the Morse theory sense, and an associated obstruction bundle of rank 1.

We have only one case left: $(\mbox{ind}_p(H_X),\mbox{ind}_q(H_X))=(1,2)$ (and both labels $\bold{hyp}$). Then $\mbox{ind}_{W_0}(u)=\mbox{ind}_M(u_0)=1$, and $u=u_0$ is non-broken, i.e.\ has only one level, with a non-trivial main level. The bottom component $u_0$ is regular inside its hypersurface, as can be checked via automatic transversality. The resulting configuration, depicted in Figure \ref{cylconf}, also has a rank 1 obstruction bundle and an evil twin, since it lies over an index $1$ flow-line connecting a hyperbolic point in $\Sigma_+$ to the minimum in $\Sigma_-$. 

One may also study the obstruction bundles for all the configurations in Figure \ref{cylconf}, and see that we are in the situation explained in Remark \ref{orbifoldcount}. We conclude that, after perturbing to generic $J$, the count of gluings of the configurations considered is zero, and finishes case 2.

\vspace{0.5cm} 

\textbf{Case 3. 1-punctured torus with one positive end.} 

\begin{figure}[t]\centering
\includegraphics[width=0.2\linewidth]{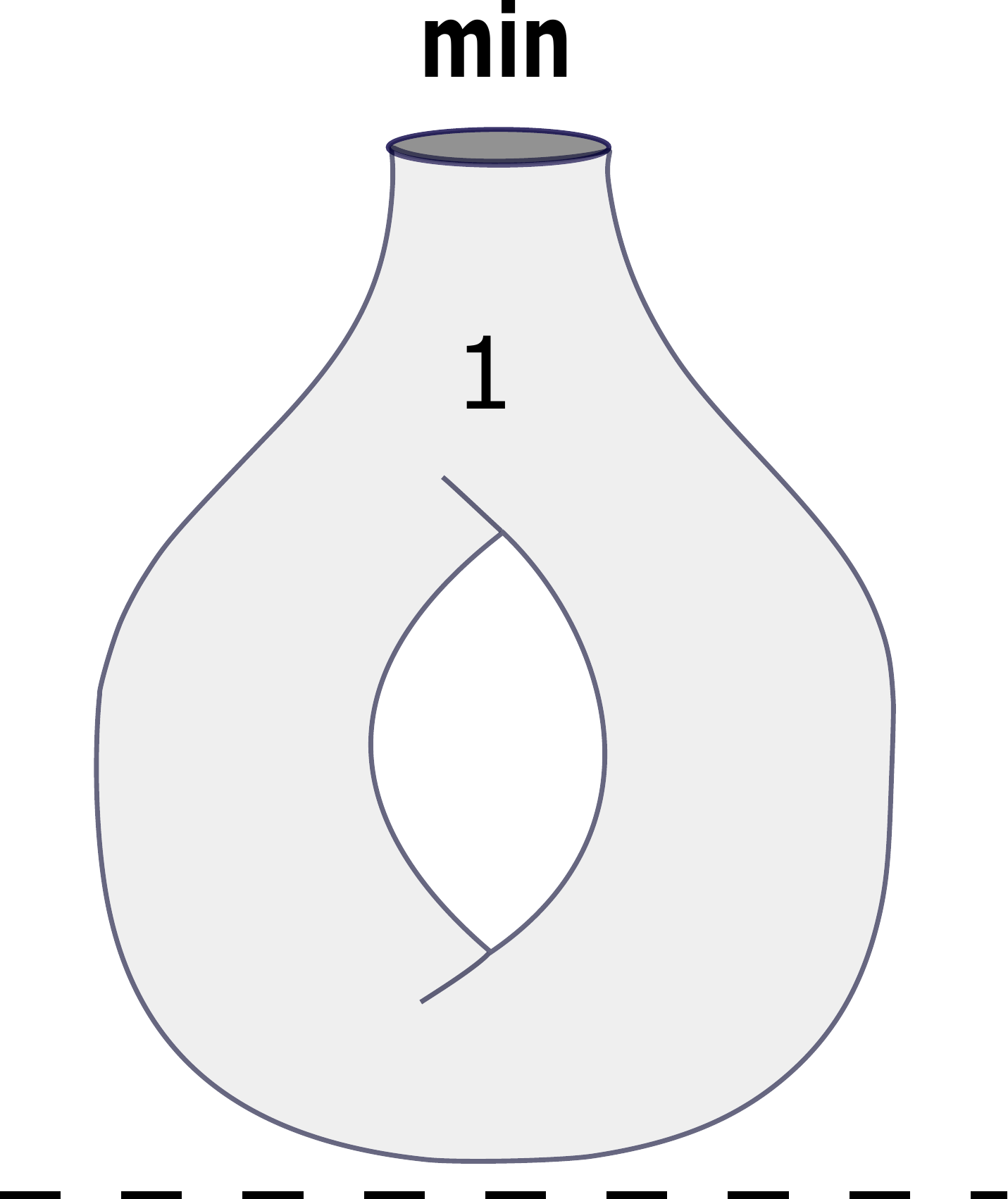}
\caption{\label{moreconfigs2} The sporadic configuration.}
\end{figure}

We will not give details of the classification, since we shall not need it for our purposes. There are 29 such configurations, all of which have obstruction bundles. Moreover, all have mirror cancelling pairs, except one of them, depicted in Figure \ref{moreconfigs2}. See \cite{Mo2}.

\subsection{The sporadic configuration}\label{esp}

In this section, we address Conjecture \ref{1tc}. Let $u$ be the sporadic configuration depicted in Figure \ref{moreconfigs2}. Let us first observe that $u$ is in fact regular in $\mathbb{R}\times M$. Indeed, since it cannot be a cover of a plane, it is somewhere injective, and so regular in the hypersurface $H:=\mathbb{R}\times H_{min}$. Equation (\ref{indN}) implies that $\mbox{ind}\;\mathbf{D}_{u;N_H}^N=0$, and since $\mathbf{D}_{u;N_H}^N$ is injective ($H$ is cylindrical), we have that it is surjective. By (\ref{cokeq}), we deduce the claim.

\subparagraph*{Some string topology.} Let us recall the definition of Goldman--Turaev string bracket and co-bracket operations, in the case of a hyperbolic surface. This is the starting point of Chas--Sullivan's string topology \cite{CS}. We follow the exposition of \cite{Sch}, which builds on \cite{CL09}.

Let $X$ be a closed surface of genus at least $2$. The set of non-trivial free homotopy classes of loops (i.e.\ a conjugacy class of the fundamental group) is countable. We choose some ordering, label by $[i]$ the $i$-th one, and make the convention that $[-i]$ denotes a change in orientation of the corresponding geodesic. Let $V$ denote the vector space generated by the free homotopy classes of non-contractible loops on $X$. We define the string co-bracket
$$
\Delta: V \rightarrow V \otimes V
$$
as follows. Represent $[i]$ by a string $s$ in general position, so that it has finite self-intersections $x_1,\dots,x_k$. At any $x_j$ there are two directions, which we orient by the orientation on $X$. We can resolve the string $s$ at each self-intersection $x_j$, obtaining two new strings $s_1^j$ and $s_2^j$. For instance, $s_1^j$ is obtained by following the first direction out of $x_j$ and taking the piece of $s$ connecting $x_j$ back to itself, and similarly for $s_2^j$. Define
$$
\Delta([i])=\sum_{j=1}^k [s_1^j]\otimes [s_2^j]-[s_2^j]\otimes [s_1^j]
$$
This co-bracket is a Lie co-bracket, i.e.\ it is bilinear, co-anti-symmetric and satisfies the co-Jacobi identity. Similarly, one define the string bracket
$$
\nabla: V\times V \rightarrow V,
$$
which is a Lie bracket (see \cite{Sch}).

In this notation, we have the coefficients $A^i_{j,k}$ for $\Delta$, defined by
$$
\Delta([i])=\sum_{j,k}A^i_{j,k}[j]\otimes [k]
$$
On the SFT side, denote by $d^i$ the SFT-count of sporadic configurations inside $\mathbb{R}\times ST^*X$, whose positive asymptotic corresponds to $[i]$. Cieliebak--Latschev show, by looking at the boundary of index 1 moduli spaces of holomorphic curves with Lagrangian boundary components in $T^*X$, the existence of linear relationships between the count of index 1 curves in $\mathbb{R}\times ST^*X$, and the coefficients of the above string operations. This allows them to compute the full SFT Hamiltonians for $T^*X$ and $\mathbb{R}\times ST^*X$. In the case of 1-punctured tori, one has the following \cite[Thm.\ 2.2.1]{Sch}: 
\begin{equation}\label{dArel}
d^i=\sum_{j} A^i_{j,-j}
\end{equation}

\subparagraph*{Relationship to 1-torsion} Figure 17 in \cite{CL09} gives an example of a closed geodesic $[i]$ for which $d^i$ is non-zero. Denote by $q_i$ the corresponding SFT generator in $\mathbb{R}\times M$. Since we have shown that these sporadic configuration is regular in the latter, we have that
$$
\mathbf{D}\left(\frac{q_i}{d^i}\right)=\hbar + \mathcal{O}(\hbar^2)
$$
Observe that we have used that there are no holomorphic disks.

Now, all the other contributions to the above differential come from potential index 1 building nodal configurations $\mathbf{u}$ with only one positive end approaching $[i]$, which glue to honest curves after introducing an abstract perturbation. If such building is entirely contained in $\mathbb{R}\times H_{min}$, then 
$$
1=\mbox{ind}_M\;(\mathbf{u})=1-\#\Gamma^-(\mathbf{u}),
$$ from which we obtain that $\#\Gamma^-(\mathbf{u})=0$. And if $\mathbf{u}$ is not contained in this cylindrical hypersurface, then from the classification of the previous section, we know that it has a twin, obtained by exploiting the symmetries of the Morse function on $\Sigma$. However, we do not know whether the latter building configurations have associated obstruction bundles, and indeed there are a priori many of them to classify (they can have any number of negative punctures, as well as any positive genus). In any case, this is strong geometric evidence for the following:

\vspace{0.5cm}

\textbf{Conjecture} Any building configuration that may contribute to the differential of $q_i$, and does not lie in $\mathbb{R}\times H_{min}$, cancels its twin when counted in SFT.

\vspace{0.5cm}

If the conjecture were true, then we have that the term $\mathcal{O}(\hbar^2)$ \emph{does not have any $q$ variables}. Since each term of the form $1+\mathcal{O}(\hbar)$ not involving $q$ variables can be inverted as a power series in $\hbar$, we obtain
$$
\mathbf{D}\left(\frac{q_i}{d^i}(1+\mathcal{O}(\hbar))^{-1}\right)=\hbar,
$$
and so our $5$-dimensional model would have 1-torsion (which would prove Conjecture \ref{1tc}).

\vspace{0.5cm}

We observe that this absence of $q$-variables is a purely $5$-dimensional phenomenon, in the following sense: Since $\mbox{ind}_{W_0}(u)=-\chi(u)$ for a curve $u$ inside $H=\mathbb{R}\times H_{min}$, the only possible contributions for the differential of $q_i$ \emph{inside} $H$, come from our sporadic configuration, as well as 3-punctured spheres with 1-positive end at $[i]$, and two negative ones, at $[j]$, and $[k]$, say. However, the latter configurations have index $-1$ in $\mathbb{R}\times M$, so do not contribute to $\mathbf{D}q_i$ in the ambient manifold. Moreover, from \cite{CL09} and \cite{Sch} we have that the SFT count of such curves inside $H$, which we denote $a^i_{j,k}$, coincides with the string coefficients $a^i_{j,k}=A^i_{j,k}$. Observe that we know from (\ref{dArel}) that there is at least one $j$ for which $A_{j,-j}^{i}$ is non-zero. So indeed we can find 3-punctured spheres in $H$ which contribute non-trivially to the SFT count in $H$, but which do not in $\mathbb{R}\times M$. 

\subsection{Non-SFT proof of Theorem \ref{nocob1}}\label{proofofnocob}

In this section, we prove Theorem \ref{nocob}. Recall the definition of a Giroux $2\pi$-torsion domain: Given a Liouville pair $(Y^{2n-1},\alpha_+,\alpha_-)$, this is the contact manifold $(GT,\xi_{GT}):=(Y\times [0,2\pi]\times S^1,\ker\lambda_{GT})$, where 
\begin{equation}
\begin{split}
\lambda_{GT}&=\frac{1+\cos(r)}{2}\alpha_++\frac{1-\cos(r)}{2}\alpha_-+\sin (r)d\theta\\
\end{split}
\end{equation} and the coordinates are $(r,\theta)\in [0,2\pi] \times S^1$.

The above contact manifold carries a suitable notion of an SOBD, and the contact form $\lambda_{GT}$ can be viewed as a Giroux form. These SOBDs, which we call \emph{Giroux SOBDs}, were described in detail in \cite{Mo,Mo2}, for more general contact manifolds obtained by gluing together collections of ``Giroux domains''. In the case of Giroux $2\pi$-torsion domains, the SOBD structure is obtained by declaring small $\delta$-collar neighbourhoods of the slices $\{r \in \{0,\pi,2\pi\}\}$ to be the ``paper'' components, and their complement, the ``spine''. The paper components are then trivial fibrations over $Y$ with fibers (the pages) which are cylinders $[-\delta,\delta]\times S^1$, and the spine components are trivial $S^1$-fibrations over a Liouville domain of the form $Y \times I$, for some interval $I$. In particular, the model A construction of \cite{Mo,Mo2} can be used. One obtains a Giroux form $\Lambda_{GT}$ which lies in the isotopy class of $\lambda_{GT}$, together with a finite energy foliation of $\mathbb{R}\times GT$, such that the cylindrical pages lift as holomorphic cylinders with two positive ends asymptotic to Reeb orbits corresponding to critical points of a Morse function in $Y \times I$. Moreover, we have a uniqueness theorem for punctured holomorphic curves in $\mathbb{R}\times Y$, which states that any other holomorphic curve, with asymptotics which are simply covered and correspond to critical points, has to be a reparametrization of a holomorphic cylinder lifting a page (see Theorem 3.9 in \cite{Mo2}, and its adaptation in the proof of Theorem 5.2 also in that thesis). If one chooses a critical point of index $1$, and a critical point of index $2n$, there exists a unique ($\mathbb{R}$-translation class of) a regular index $1$ holomorphic cylinder in the finite energy foliation with asymptotics corresponding to these critical points, which we call $u_0$. See \cite{Mo2} for the full details.  

\vspace{0.5cm}

We are now ready to obtain an application from our knowledge of the SFT algebra of $ST^*X\times \Sigma$. Recall that we distinguish between model A (constructed in \cite{Mo2}), and model B (constructed in this paper), which give the same isotopic contact structures.

\begin{proof}[Thm. \ref{nocob1}] 

Let $(M_0,\xi_0)$ be a $5$-dimensional contact manifold with Giroux torsion, and let $i:(GT,\xi_{GT})\hookrightarrow (M_0,\xi_0)$ be a contact embedding. Consider the $5$-dimensional model $(M=ST^*X\times \Sigma,\xi)$, with $k\geq 3$, where we view $\xi=\ker \Lambda_\epsilon$ as the contact structure induced by a model B contact form $\Lambda_\epsilon$. The $\epsilon$-parameter does not change the isotopy class, and we will choose it suitably below. We may take a contact form $\Lambda_0$ for $\xi_0$, such that it coincides with a model A contact form $\Lambda_{GT}$ on $i(GT)$ as described above. Assume that $(W,d\lambda)$ is an exact cobordism with convex end $(M_0,\xi_0)$ and concave end $(M,\xi)$. Let $J_1$ be a model B almost complex structure on the negative half-symplectization of $(M,\xi)$ (which is compatible with $\Lambda_\epsilon$). Let also $J_0$ be a $\Lambda_0$-compatible almost complex structure on the positive half-symplectization of $(M_0,\xi_0)$, which coincides with a model A almost complex structure on $i(GT)$ for which we have the finite energy foliation of $\mathbb{R}\times i(GT)$ described before. 

By attaching small trivial symplectic cobordisms to $W$, we can assume that $\lambda\vert_{M_0}$ and $\lambda\vert_{M}$ are positive and large constant multiples of $\Lambda_0$ and $\Lambda_\epsilon$, respectively, where the constant can be chosen $\epsilon$-independent. Consider the Liouville completion $\widehat{W}$, obtained from $W$ by attaching one positive and one negative cylindrical ends, and the natural extension of $\lambda$ to $\widehat{W}$. Choose a $\lambda$-compatible almost complex structure $J$, such that $J$ coincides with $J_0$ and $J_1$ on their respective ends, and such that $J$ is generic along the main level. Consider also the compactification $\overline{W}^\infty$, obtained by adding $+\infty$ and $-\infty$ to $\widehat{W}$, on which we have extensions of both the symplectic form $d\lambda$ and $J$ (see \cite{Wen4} for more details on this construction). 

A we have already described, there exists a distinguished ($\mathbb{R}$-translation class of a) regular index 1 $J_0$-holomorphic cylinder $u_0$, lying in the upper levels of $\widehat{W}$. Let $T$ be the total action of its two positive ends. Recall that by the model B construction, there is an action threshold $T_\epsilon>0$, such that $\lim_{\epsilon\rightarrow 0}T_\epsilon=+\infty$ and such that every $\Lambda_\epsilon$-closed orbit with action less than $T_\epsilon$ corresponds to a critical point in $M$, or lies along the spine $Y \times \mathcal{U}$. Observe that this construction also holds if we multiply the contact form by a constant positive number which is $\epsilon$-independent. Take $\epsilon>0$ small enough so that $T_\epsilon>T$. By Lemma \ref{sympcob}, by having attached a suitable exact symplectic cobordism to the concave end of $W$ prior to taking the completion, we can assume that $\lambda\vert_{M}$ is a constant ($\epsilon$-independent) multiple of the model B contact form $\Lambda_\epsilon$ corresponding to this particular $\epsilon>0$. By Stokes' theorem, we get that the action of the Reeb orbits appearing in any building configuration arising as a limit of a curve in the moduli space of $u$ is bounded by $T_\epsilon$. Therefore, in what follows, we can make use of the discussion in Section \ref{proooof}, from which we gather that, whenever this cylinder breaks in lower levels, the Reeb orbits we obtain all correspond to critical points, and none lies over the spine.

Denote by $\overline{\mathcal{M}}$ the connected component containing $u_0$ of the nodal compactification of the space of finite energy $J$-holomorphic curves in $\overline{W}^\infty$. It is a 1-dimensional compact manifold with boundary, having one boundary component corresponding to the projection of $u_0$ to $M_0$.

We observe first that no element of $\overline{\mathcal{M}}$ can break in an upper level. Indeed, by the aforementioned uniqueness Theorem 3.9 in \cite{Mo2} (more precisely, its adaptation as in the proof of Theorem 5.2 in \cite{Mo2}) all of the upper components of a stable building $u \in \partial\overline{\mathcal{M}}$ lying in $\mathbb{R}\times M_0$ have to correspond to flow-line cylinders of a Morse function (which have non-negative index). But then $u$ can only consist of an index $1$ cylinder in the bottom end, and a chain of index $0$ cylinders on top of it, which are necessarily trivial. Since $u$ is stable, this configuration is not allowed unless $u$ coincides with its bottom level, but then $u$ is non-broken.  

Similarly, if an element in $\overline{\mathcal{M}}$ breaks with a non-trivial main level component, and no lower level ones, then uniqueness gives that there is only one possible breaking configuration. It consists of a single upper level, having an index $1$ somewhere injective flow line cylinder, together with a trivial (index zero) cylinder; and a nontrivial somewhere injective index $0$ main level component consisting of a cylinder with two positive ends, glued along simply covered Reeb orbits corresponding to critical points of Morse index $1$ and $3$. This configuration has an evil twin, since the index $1$ component does (this is analogous to what we have already observed for model B), and it also consists of somewhere injective components. Therefore, we may glue this configuration to a unique honest cylinder. Then we obtain two $1$-dimensional moduli spaces, each having an open end on each twin configuration. We may then canonically identify them along these open ends, and obtain a new moduli space, which we still call $\overline{\mathcal{M}}$, for which we repeat this process.

Assume that an element in $\overline{\mathcal{M}}$ breaks in a lower level of $\overline{W}^\infty$, corresponding to the negative symplectization of $M$. Since the cobordism is exact, there are no holomorphic caps, so that all components of all possible breaking configurations are still cylinders. Then we can argue similarly as above, but using an obstruction bundle. From the discussion in Section \ref{proooof}, we get an obstruction bundle and a twin. Recall also that we have also made sure that we are in the situation of Remark \ref{orbifoldcount}. We may choose finite and sufficiently large gluing parameters, such that the twin configuration glues to honest holomorphic cylinders, whose number is the algebraic count of zeroes of the section of the corresponding obstruction bundle. This number, which is independent of the gluing parameters as long as they are large enough, is the opposite of the original configuration. Since we know by construction that the latter glues at least once, its evil twin will also. And we can proceed as before. Observe that we haven't needed an abstract perturbation.

After these identifications, we will have constructed a 1-dimensional moduli space which has only one boundary component, which is a contradiction. This finishes the proof.
\end{proof} 

\appendix

\section{Super-rigidity for punctured holomorphic curves}\label{obpants}

In this appendix, we will derive Theorem \ref{superig}, a general result bearing some independent interest. While we shall not use it in this paper, it will come as a byproduct of Proposition \ref{unb} below, which is used in the main text.

For a holomorphic degree $d$ branched cover $\varphi: (\dot{S},j)\rightarrow (\dot{\Sigma},i)$ between punctured Riemann surfaces $\dot{S}=S\backslash\Gamma(S)$ and $\dot{\Sigma}=\Sigma\backslash\Gamma(\Sigma)$, define by $$\kappa(\varphi)=\sum_{w\in\Gamma(S)}\kappa_w\geq \#\Gamma(S)$$ the \emph{total multiplicity} of $\varphi$, where $\kappa_w=Z(d_w\varphi)+1$ is the multiplicity of $\varphi$ at the puncture $w$ (where $Z(d_w\varphi)$ is the vanishing order of $d\varphi$ at $w$). Observe that 
$$
\kappa(\varphi)=\sum_{z \in \Gamma(\Sigma)}\left(\sum_{w \in \varphi^{-1}(z)}\kappa_z\right)=d\#\Gamma(\Sigma)
$$  
Therefore $d\#\Gamma(\Sigma)\geq \#\Gamma(S)$, with equality if and only if there is no branching at the punctures. We will refer to the sum of the branching of $\varphi$ at interior points of $\dot{S}$ with the branching at the punctures as the \emph{total branching} of $\varphi$.

\begin{prop}\label{unb} Let $v:(\dot{\Sigma},i)\rightarrow (W,J)$ be a somewhere injective possibly punctured holomorphic curve in a 4-dimensional symplectic cobordism, with index zero and $\chi(v)=0$. Assume that there exist a trivialization $\tau$ for which the asympotics of $v$ have vanishing Conley--Zehnder index. Assume also that $v$ satisfies $\ker \mathbf{D}_v^N = 0$.  Then $\mathbf{D}_u^N$ is also injective for any multiple cover $u=v \circ \varphi$ of $v$, such that $\varphi$ has strictly positive total branching. 
\end{prop}

The proof follows by an adaptation of Hutchings' \emph{magic trick} (see \cite{Hut2}), as used in Prop.\ 7.2 of \cite{Wen6} in the case of closed curves. We redo Wendl's computation, and the unpunctured case recovers Wendl's results for the torus.

In the proof of the proposition, we will use Siefring's intersection pairing, and the adjunction formula for punctured curves. We refer to \cite{Wen8} (especially Theorem 4.4) for the statement, and details on the necessary definitions. In particular, the symbol $\overline{\sigma}(u)$ denotes the \emph{spectral covering number} of $u$. This is defined as the sum, over all punctures, of the covering multiplicity of any eigenfunction, corresponding to the eigenspaces of the extremal eigenvalues, of the associated asymptotic operators. Since the multiplicity of the eigenfunction divides that of the asymptotic orbit, we have $\overline{\sigma}(u)=\#\Gamma(u)$ whenever the asymptotics of $u$ are simply covered. Also, the symbol $\delta(u)$ denotes the algebraic count of double points and critical points of $u$. It vanishes if and only if $u$ is embedded. The symbol $\delta_\infty(u)$ is the algebraic count of ``hidden double points at infinity''. If different punctures of $u$ asymptote different orbits, and all of them are simply covered, then this number vanishes (see Thm. 4.17 in \cite{Wen 8} for the general case).

\begin{proof}(Prop. \ref{unb}) As explained e.g.\ in App.\ of \cite{Wen6}, one can construct an almost complex structure $J_N$ on the total space of the normal bundle $\pi: N_v \rightarrow \dot{\Sigma}$, such that there is a 1-1 correspondence between holomorphic curves $u_{\eta}:(\dot{S},j) \rightarrow (N_v,J_N)$ and sections $\eta \in \ker \mathbf{D}_u$ along holomorphic branched covers $\varphi=\pi \circ u_\eta: (\dot{S},j)\rightarrow (\dot{\Sigma},i)$. 

Let $u=v\circ\varphi$ be a branched cover of $v$. We will assume nothing on $\chi(v)$ until the very end. Let $d:=\mbox{deg}(\varphi)\geq 1$. We proceed by induction on $d$. The base case $d=1$ holds by assumption, and so take $d > 1$. Assume by contradiction that there exists a nonzero $\eta \in \ker \mathbf{D}_u^N$, and take the corresponding $u_\eta$. We can view $v$ as a $J_N$-holomorphic embedding inside $N_v$ (as the zero section), and we have that $u_\eta$ is homologous to $d[v]$, where $[v]$ is a \emph{relative} homology class. By the induction hypothesis, we can assume that $u_\eta$ is somewhere injective. Observe that $v$, as a holomorphic map to $N_v$, is embedded, and each of its asymptotics is distinct and simply covered, so that we get $\delta(v)=\delta_\infty(v)=0$, and $\overline{\sigma}(v)=\#\Gamma(v)$.

Since the asymptotics of $v$ have vanishing Conley--Zehnder index, their extremal winding numbers also vanish. Therefore the extremal eigenfunctions are constant, so that their covering multiplicity coincides with that of the corresponding asymptotic orbit. If we denote by $\kappa(u_\eta)$ the sum of all the multiplicities of the asymptotics of $u_\eta$ (which coincides with $\kappa(u):=\kappa(\varphi)$), we obtain that the spectral covering number $\overline{\sigma}(u_\eta)=\kappa(u_\eta)=\kappa(u)=d\#\Gamma(v)$.

We can then compute its Siefring self-intersection number:
\begin{equation}
\label{uno}
\begin{split}
u_\eta*u_\eta&=2(\delta(u_\eta)+\delta_\infty(u_\eta))+c_1^\tau(u_\eta^*TN_v)-\chi(u_\eta)+\overline{\sigma}(u_\eta)-\#\Gamma(u_\eta)\\&=2(\delta(u_\eta)+\delta_\infty(u_\eta))+dc_1^\tau(v^*TN_v)-\chi(u_\eta)+d\#\Gamma(v)-\#\Gamma(u)
\\
\end{split}
\end{equation}

Moreover, Riemann--Hurwitz yields
\begin{equation}\label{cinco}
\chi(u_\eta)=\chi(u)=d\chi(v)-Z(d\varphi)
\end{equation}
Identifying the normal bundle to $v$ inside $N_v$ with $N_v$ itself, the adjunction formula then produces
\begin{equation}\label{tres}
c_1^\tau(v^*TN_v)=\chi(v)+c_1^\tau(N_v)
\end{equation}
Moreover, we have that
\begin{equation}\label{indzeroN}
0=\mbox{ind}\; \mathbf{D}_v^N = \chi(v) + 2c_1^\tau(N_v)
\end{equation}
Using (\ref{indzeroN}) and the vanishing of the Conley-Zehnder indices, we obtain 
\begin{equation}\label{dos}
\begin{split}
u_\eta*u_\eta&=d^2v*v\\
&=d^2(c_1^\tau(v^*TN_v)-\chi(v))\\
&=d^2 c_1^\tau(N_v)\\
&=-\frac{d^2}{2}\chi(v)
\end{split}
\end{equation}

Combining (\ref{uno}), (\ref{cinco}), (\ref{tres}), (\ref{indzeroN}) and (\ref{dos}), we get
\begin{equation}
\begin{split}
0\leq & 2(\delta(u_\eta)+\delta_\infty(u_\delta))\\
=&-\frac{d^2}{2}\chi(v)-d(\chi(v)+c_1^\tau(N_v))+d\chi(v)-Z(d\varphi)+\#\Gamma(u)-d\#\Gamma(v)\\
=&\frac{d(1-d)}{2}\chi(v)-Z(d\varphi)+\#\Gamma(u)-d\#\Gamma(v)
\end{split}
\end{equation}
If we assume that $\chi(v)=0$ (or $\chi(v)\geq 0$), the above is negative, and we obtain a contradiction. 
\end{proof}

In order to obtain Theorem \ref{superig}, we need another result. Recall that a closed hyperbolic Reeb orbit is one whose linearized return map does not have eigenvalues in the unit circle. In dimension $3$, if their Conley--Zehnder index is even (in one and hence every trivialization), then they are hyperbolic.

\begin{thm}\cite{Wen7}\label{PM} For generic $J$, all index $0$ punctured holomorphic curves in a $4$-dimensional symplectic cobordism, which are unbranched multiple covers, which have zero Euler charactersitic, and such that all of their asymptotics are hyperbolic, are Fredholm regular. 
\end{thm}

This is proved by the same arguments as in \cite{Wen6} (which only deals with closed curves), and we will not include a proof. Although this paper has recently been withdrawn due to a gap in a lemma, which was crucial for the super-rigidity results there claimed, the hypothesis included in Theorem \ref{PM} bypasses that gap. Indeed, the problem in the paper has to do with the failure in general of \emph{Petri's condition} (see the recent blog post \cite{Wen9} explaining this). In the current situation, for a curve $u$ as in Theorem \ref{PM}, the normal Cauchy-Riemann operator $\mathbf{D}_u^N$ is defined on a line bundle $N_u$ whose adjusted first Chern number is non-positive: we have $2-2g-\#\Gamma=0$, so that $2-2g-\#\Gamma_{odd}=\#\Gamma_{even}$, and therefore
$$
2c_1(N_u,u)=\mbox{ind}(u)-2+2g+\#\Gamma_{even}=-\#\Gamma_{odd}\leq 0 
$$  
This means that either $\ker \mathbf{D}_u^N=0$, or every element in $\ker \mathbf{D}_u^N$ is nowhere vanishing. Therefore Petri's condition is vacuously satisfied (this is also pointed out in the aforementioned blog post). 

With this said, Theorem D in that paper needs to be generalized to the setting of asymptotically cylindrical curves in cobordisms. This is mostly a matter of putting it in the proper functional-analytic setup for punctured curves, and almost nothing else changes. Since all orbits are hyperbolic, Fredholm indices of multiple covers are related to indices of the underlying somewhere injective curves via the same multiplicative relations (involving the Riemann--Hurwitz formula) as in the closed case. In fact, the multiplicativity of the Conley--Zehnder indices is all what is needed. With that in mind, the same dimension counting argument as in Theorem B of said paper works in the case of hyperbolic punctures.

\begin{proof} (Super-rigidity, Theorem \ref{superig}) In the case where $v$ is as in the hypothesis of Theorem \ref{superig}, generically we can assume it is immersed and regular. Since $\mbox{ind}\;\mathbf{D}_v^N=0$, then $\ker \mathbf{D}_v^N=0$ by regularity, and so the hypothesis of Proposition \ref{unb} are satisfied for generic $J$. If we also assume that its totally unbranched covers satisfy the conclusion of the proposition (which is also a generic condition by Theorem \ref{PM}), the result follows.
\end{proof}

\section{Lutz--Mori twists and Giroux torsion}
\label{LMtwists}

This appendix is devoted to describe the Lutz--Mori twists as defined in \cite{MNW}, and to use them to model A of \cite{Mo}, as described in the introduction. In our particular case, one can use CCH (cylindrical contact homology) to distinguish the twisted model A contact structures for different amounts of twisting, something which has been done already in \cite{MNW} (Theorem 8.13) for very similar models, so we will not include any details.  

\subsection{The twisted contact structures}

Reccall that a hypersurface $H\subseteq M$ in a contact manifold $(M,\xi)$ is a \emph{$\xi$-round hypersurface modeled on $(Y,\xi_0)$} if it is transverse to $\xi$ and admits an orientation preserving identification with $S^1\times Y$, for some contact manifold $(Y,\xi_0)$, such that $\xi \cap TH = TS^1 \oplus \xi_0$. This notion is defined in \cite{MNW}.

Take $(Y^{2n-1},\alpha_-,\alpha_+)$ a Liouville pair, and denote by $\xi_\pm =\ker \alpha_\pm$. Consider the \emph{Giroux $2\pi l$-torsion domain} modeled on $(Y,\alpha_-,\alpha_+)$, which is the contact manifold $$GT_l^+:=(Y \times [0,2\pi l] \times S^1,\xi_{GT}=\ker\lambda_{GT}),$$ where $$\lambda_{GT}=\frac{1+\cos (r)}{2}\alpha_++\frac{1-\cos (r)}{2}\alpha_-+\sin (r) d\theta,$$ both of whose boundary components are $\xi_{GT}$-round hypersurfaces modeled on $(Y,\xi_+)$.  We can also consider the modified version given by $$GT_l^-:=(Y \times [\pi/2,2\pi l+\pi/2] \times S^1,\xi_{GT}=\ker\lambda_{GT}),$$ whose boundaries are now modeled on $(Y,\xi_-)$. 

\vspace{0.5cm}

For the model A version of the contact structure $\xi_k$ in $M=Y\times \Sigma$ (in the same isotopy class as our model B contact structures), we have that $\xi_k=\xi_\pm \oplus T\Sigma_\pm$ over the regions $Y \times \Sigma_\pm$. Therefore, we obtain ($k$ copies of) round hypersurfaces $H_\pm$ modeled on $(Y,\xi_\pm)$, corresponding respectively to the boundary components of $Y \times \Sigma_\pm$. We may therefore perform an $l$-fold \emph{Lutz--Mori twist} along the $H_\pm$. This can be done in two equivalent ways. The first consists in cutting our model along each of the $k$ copies of $H_-$ and gluing a copy of $GT_l^-$, and the second in doing the same for $H_+$, along which we glue $GT_l^+$. This yields a manifold which is diffeomorphic to $M$, with a contact structure on $M$ which we denote $\xi^{l}_k$. This contact structure is homotopic as almost contact structures to the original contact structure $\xi_k$. Since $\xi^{l}_k$ has Giroux torsion, Theorem \ref{AlgvsGiroux} implies that it has algebraic $1$-torsion (for $l>0$), so that they are indistinguishable from this invariant alone. Corollary \ref{nocob} follows. On the other hand, one may use cylindrical contact homology to distinguish them \cite{Mo2}.   

\newpage

\bibliographystyle{alpha}

\end{document}